\documentclass{amsart}     

\usepackage{amssymb}
\usepackage{graphicx}
\usepackage{amscd}
\usepackage{mathptmx}
\newtheorem{thm}{Theorem}[section]
\newtheorem{lm}[thm]{Lemma}
\newtheorem{prop}[thm]{Proposition}
\newtheorem{crl}[thm]{Corollary}

\theoremstyle{definition}
\newtheorem{df}[thm]{Definition}
\newtheorem{ex}[thm]{Example}
\theoremstyle{remark}
\newtheorem{rem}[thm]{Remark}

\newcommand{\coker}{\operatorname{coker}}

\newcommand{\wt}{\operatorname{wt}}
\newcommand{\diag}{\operatorname{diag}}
\newcommand{\irightarrow}{\stackrel{\sim}{\longrightarrow}}

\newcommand{\M}{{\mathcal M}}

\newcommand{\C}{{\mathbb {C}}}

\newcommand{\Z}{{\mathbb {Z}}}
\newcommand{\R}{{\mathbb {R}}}

\newcommand{\D}{{\mathcal D}}

\newcommand{\z}{{\mathbf z}}

\newcommand{\m}{{\mathbf m}}

\renewcommand{\O}{{\mathcal{O}}}

\newcommand{\tu}{{\tilde{u}}}

\newcommand{\tf}{{\tilde{f}}}

\newcommand{\ind}{\operatorname{index}}

\begin{document}

\title[Spin Equations]{Geometry and analysis of spin equations}

\author{Huijun Fan}
\address{Peking University and Max-Planck Institute for Mathematics, Leipzig}
\thanks{The first author was partially supported by the Research Fund for returned overseas Chinese Scholars 20010107 and the Partner Group of the Max-Planck
Institute for mathematics in the sciences and the Chinese Academy of
Sciences, and later by NSFC 10401001, NSFC 10321001 and NSFC
10631050. The second author was partially funded by NSF Grant
DMS0105788.  The third author was partially supported by an NSF
grant and a Hong Kong RGC grant.}
\author{Tyler J. Jarvis}
\address{Brigham Young University}
\author{Yongbin Ruan}
\address{University of Wisconsin-Madison}

\begin{abstract}
We introduce $W$-spin structures on a Riemann surface $\Sigma$ and
give a precise definition to the corresponding $W$-spin equations
for any quasi-homogeneous polynomial $W$. Then, we construct
examples of nonzero solutions of spin equations in the presence of
Ramond marked points. The main result of the paper is a compactness
theorem for the moduli space of the solutions of $W$-spin equations
when $W=W(x_1,\dots,x_t)$ is a non-degenerate quasi-homogeneous
polynomial with fractional degrees (or weights) $q_i<1/2$ for all
$i$. In particular,  the compactness theorem holds for the
superpotentials $E_6,E_7,E_8$, or $A_{n-1}, D_{n+1}$ for $n\ge 3$.
\end{abstract}

\maketitle

\tableofcontents

\section{Introduction}

Since Donaldson's remarkable work in the 1980's regarding anti-self
dual equations on 4-manifolds, the study of the solution space
(moduli space) of a nonlinear elliptic PDE has attracted a great
deal of attention in geometry and topology. For example, the famous
Donaldson invariants, Seiberg-Witten invariants and Gromov-Witten
invariants were constructed out of the moduli space of anti-self
dual Yang-Mill's equations, the Seiberg-Witten equations, and the
$\bar{\partial}$ equation, respectively. These invariants have
revolutionized many subjects of mathematics. There is a well-known
program originated by  Donaldson and Taubes to construct invariants
out of an elliptic equation. The first, and often most important,
step is to study the analytic properties of equations such as
compactness and Fredholm theory. Once the elliptic equation under
study possesses nice analytic properties, we can apply sophisticated
machinery to extract topological invariants.  It is often a
difficult problem to show that an elliptic equation has ``nice''
analytic properties.
 In terms of the analytic framework, the best
understood ones are conformally invariant first order equations such
as the anti-self dual equation on a 4-manifold and the
$\bar{\partial}$-equation in two dimensions. One can understand
these two types of equations through Uhlenbeck's bubbling analysis.
Beyond conformally invariant equations, there is not yet a standard
framework. In this paper, we study  a new type of elliptic equation
(spin equation) on Riemann surfaces introduced by Witten. A spin
equation is of the form
$$\bar{\partial}u_i+\overline{\frac{\partial W}{\partial
u_i}}=0,$$
    where $W$ is a quasi-homogeneous polynomial, and $u_i$ is
    interpreted as the section of an appropriate orbifold line bundle
    on a Riemann surface $\Sigma$. Typically, a spin equation is not
    conformally invariant.

    The simplest spin equation ($A_{r-1}$ case) is of the form
    $$\bar{\partial}u+r\bar{u}^{r-1}=0.$$
    It was introduced by
    Witten \cite{W2} more than ten years ago as a generalization
    of topological gravity. Somehow, it was buried in the
    literature without attracting much attention.
    Recently, Witten generalized it to
    a spin equation for an arbitrary quasi-homogeneous polynomial \cite{W3}
    and coined it the ``Landau-Ginzburg A-model.'' There appears to be a mirror symmetry between this theory
    and usual Landau-Ginzburg theory (B-model) \cite{Ka1,Ka2}. The construction of invariants and the application to
    mirror symmetry will be left to a separate article. We will focus on the analytic
    aspect of the theory.

    We begin with a brief history of the motivation behind the spin
    equation.
    Around 1990, Witten proposed a remarkable conjecture relating
    the intersection number of the Deligne-Mumford moduli space of
    stable curves with the KdV integrable hierarchy \cite{W1}. His conjecture
    was soon proved by Kontsevich \cite{K}. About the same time, Witten also proposed a
    generalization of his conjecture. In his generalization, the stable curve
    is replaced by a certain root of the canonical bundle (spin-curve) and the KdV-hierarchy
    was replaced by certain, more general, KP-hierarchies called nKdV, or Gelfand-Dikii, hierarchies. Since then, the moduli
    space of spin-curves has been rigorously constructed by the second author and his collaborators
    \cite{AJ,J1,J2,JKV}.

    An important phenomenon in the theory of spin curves is the appearance of
    Neveu-Schwarz and Ramond marked points. Recall that $u_i$ is
    the section of a certain orbifold line bundle $L_i$. Assume that
    all the orbifold points are marked points. A marked point with
    trivial orbifold structure is called a Ramond marked point, and
    otherwise it is called a Neveu-Schwarz marked point. Contrary to intuition,
    Ramond marked points are much harder to study than Neveu-Schwarz
    marked points. If there is no Ramond marked point, a simple
    lemma of Witten's shows that the spin equation has only the zero
    solution. Therefore, our moduli problem becomes an algebraic
    geometry problem. In fact, Witten conjectured that
    the contribution of a Ramond marked point to the corresponding field theory is zero in the $A_{r-1}$ case (the decoupling
    of the Ramond sector).  This was proved true for genus zero in \cite{JKV} and for higher genus in \cite{P}.
    This is partly why the moduli space of spin curves has been around
    for a long time while the spin equation seems to have been lost in the
    literature. In the course of our investigation, we discovered that in the $D_{n}$-case the Ramond
    sector gives a nonzero contribution. Hence, we have to develop
    a theory to account for the contribution of the solution of the spin
    equation in the presence of Ramond marked points. The situation
    would be relatively simple if one could show that the solution has to be zero
    as in the case of pure Neveu-Schwarz marked points. Unfortunately, this is not the case, as
    we were able to construct a nonzero
    solution of the spin equation. Therefore, we have no choice but
    to develop a full moduli theory of spin equations with Ramond
    marked points. This is the first of a series of articles to
    accomplish that task.

Define $\M_r(\Sigma,W)$ to be the space of regular solutions
$(u_1,\dots,u_t)$ of the $W$-spin equations on $\Sigma$, and let
$\M(\Sigma,W)$ be the solution space of $W$-spin equations which
contains both the regular and the singular solutions. Now we can
formulate our main theorems as follows:

\begin{thm}[\bf Inner compactness]\label{main-thm1} Suppose $W=W(x_1,\dots,x_t)$ to be
a non-degenerate quasi-homogeneous polynomial with fractional
degrees (or weights $wt(x_i)$) $q_i=wt(x_i)\le 1/2$. Let
$E(u):=\sum_{z_l:\text{Ramond}}\sum_{j:W_j \text{is Ramond at}\,
z_l}Res_{z_l}(W_j(u_1,\dots,u_t))$ be the residue map from
$\M_r(\Sigma, W)$ to ${\mathbb C}$. Then
\begin{enumerate}
\item[(1)] for any $a\in {\mathbb C}\setminus [0,\infty)$, we have $E^{-1}(a)=\emptyset$; and
\item[(2)] for any $a\in [0,\infty)$, $E^{-1}([0,a])$ is a compact
space in the $L^p_1$ topology for $2\le p<\frac{2}{1-\delta}$, where
$\delta=\min\{q_1,\dots,q_t\}$.
\end{enumerate}
\end{thm}

\begin{thm}[\bf Weak compactness]\label{main-thm2} Suppose $W=W(x_1,\dots,x_t)$ to be
a non-degenerate quasi-homogeneous polynomial with fractional
degrees (or weights $wt(x_i)$) $q_i=wt(x_i)<1/2$. For any
$\epsilon>0$, let $\Sigma_\epsilon$ be the subset of $\Sigma$
consisting of points which are distance at least $\epsilon$ away
from Ramond marked points. Then the restriction of $\M(\Sigma, W)$
to $\Sigma_\epsilon$ is compact in $L^p_1(\Sigma_\epsilon)$ for
$2\le p<\frac{2}{1-\delta}$, where $\delta=\min\{q_1,\dots,q_t\}$.
\end{thm}

\begin{thm}[\bf Strong compactness]\label{main-thm3}Suppose $W=W(x_1,\dots,x_t)$ to be
a non-degenerate quasi-homogeneous polynomial with fractional
degrees (or weights $wt(x_i)$) $q_i=wt(x_i)<1/2$. Then
$\M(\Sigma,W)$ is compact with respect to the topology in
$\substack{\underbrace{L^{p_1}\times\cdots\times L^{p_t}}\\ t}$ for
$0<p_i<2(1-2q_i)/q_i$. In particular, we have
\begin{enumerate}
\item If $W$ is the $A_{r-1}$-superpotential, then $\M(\Sigma,W)$ is
compact with respect to the topology in $L^p$ for $1<p<2(r-2)$.
\item If $W$ is the $D_n$-superpotential for $n\ge 4$, then
$\M(\Sigma,W)$ is compact with respect to the topology in
$L^{p_1}\times L^{p_2}$ for $1<p_1<2(n-2), 0<p_2<4/(n-1)$.
\item If $W$ is the $E_6$ superpotential, then $\M(\Sigma,W)$ is
compact with respect to the topology in $L^{p_1}\times L^{p_2}$ for
$1<p_1<2, 1<p_2<4$.
\item If $W$
is the $E_7$ superpotential, then $\M(\Sigma,W)$ is compact with
respect to the topology in $L^{p_1}\times L^{p_2}$ for $1<p_1<2,
1<p_2<5$.
\item If $W$ is the $E_8$ superpotential, then
$\M(\Sigma,W)$ is compact with respect to the topology in
$L^{p_1}\times L^{p_2}$ for $1<p_1<2, 1<p_2<6$.
\end{enumerate}
\end{thm}

\begin{rem}
Because of the $W$-spin structure, $W(u_1,\dots,u_t)$ takes value in
the log-canonical bundle, hence we can take the residue at the
Ramond marked points. Theorem \ref{main-thm1} means that a sequence
of points $\{u^l=(u^l_1,\dots, u^l_t)\}$ in $\M_r(\Sigma, W)$ loses
compactness if and only if $E(u^l)\rightarrow +\infty$. The crucial
identity connecting the total residue and the natural energy of the
sections will be given in Theorem \ref{thm-resi-ener}. Theorem
\ref{main-thm2} means that all the possible blow-up points of a
sequence of solutions in $\M(\Sigma, W)$ are Ramond marked points.
Theorem \ref{main-thm3} means that we can compactify the space
$\M_r(\Sigma, W)$ in some $L^p$ topology so that its closure is just
$\M(\Sigma, W)$.\end{rem}

\begin{rem} Our requirement for the fractional degrees of
$W$ is sharp. If the fractional degree $q_i=1/2$ for some $i$, then
the weak and the strong compactness theorems may not be true. For
example, we can consider the $A_1$-case, where the fractional degree
is $1/2$ and the spin equation is
$$\bar{\partial}u+r\bar{u}=0.$$
This is a real-linear equation: if $u$ is a solution, then $\lambda
u$ is also a solution for any $\lambda\in {\mathbb R}$. Thus Theorem
\ref{main-thm2} and Theorem \ref{main-thm3} do not hold.
\end{rem}

    This paper is organized as follows. In Section 2, we will
    introduce the $W$-spin equations and their basic properties. In Section 3, we
    will present our key example of a nonzero solution of the spin equation
    in the presence of two Ramond marked points. More importantly, we
    construct a sequence of regular nonzero solutions converging to
    a singular one. This indicates the complexity of the problem.
 In Section 4,
we give the $L^p$ estimate of the $\bar{\partial}$ operator in
certain natural weighted Sobolev spaces arising in the problem and
prove $\bar{\partial}$ is a Fredholm operator under some mild
constraints. The core of the paper is in the last two sections,
where we establish the three compactness theorems. We will prove
Theorem \ref{main-thm1} in Section \ref{inne-comp} and Theorem
\ref{main-thm2}--\ref{main-thm3} in Section \ref{stro-comp}.

\section{Spin structures on orbicurves and spin equations}
In this section, we will introduce $W$-spin structures on
orbicurves, where $W\in {\mathbb C}[x_1,\dots,x_t]$ is a
non-degenerate  quasi-homogeneous polynomial. By means of $W$-spin
structures, one can define the $W$-spin equations on orbicurves.

Let $W\in {\mathbb C}[x_1,\dots,x_t]$ be a quasi-homogeneous
polynomial, i.e., there exist degrees $d, k_1,\dots,k_t\in \Z^{>0}$
such that for any $\lambda\in {\mathbb C}^*$
$$
W(\lambda^{k_1}x_1,\dots,\lambda^{k_t}x_t)=\lambda^d
W(x_1,\dots,x_t).
$$

\begin{df} $W$ is called \emph{nondegenerate} if
\begin{item}
\item{(1)} the fractional degrees $q_i=\frac{k_i}{d}$ are uniquely
determined by $W$; and

\item{(2)} the hypersurface defined by $W$ in weighted projective
space is non-singular, or equivalently, the affine hypersurface
defined by $W$ has an isolated singularity at the origin.
\end{item}
\end{df}

\begin{rem}
Although the first condition does not play an essential role in this
paper, it is essential for producing a compact moduli space of
orbicurves with spin structure, which is a major part of the
motivation for this paper.  Thus we will not hesitate to use this
assumption whenever it is convenient.
\end{rem}

From now on, we always assume the quasi-homogeneous polynomial $W$
is nondegenerate and the corresponding degrees $d,k_1,\dots,k_t$ of
$W$ are the {\em least} positive integer degrees.

\begin{ex}

\

\begin{enumerate}
\item We call $W(x)=x^r$ the $A_{r-1}$ superpotential.   In this case, the variable $x$ has weight $q_x=1/r$.
\item We call $W(x,y)=x^n+x y^2$ the $D_{n+1}$ superpotential; it has $q_x=1/n$ and $q_y=\frac{n-1}{2n}$.
\item Similarly, we call
$W(x,y)=x^3+y^4, W(x,y)= x^3+x y^3,$ and $W(x,y)= x^3+y^5$,  the
$E_6, E_7, $ and $E_8$ superpotentials, respectively, and their
weights are $(q_x,q_y)=(1/3,1/4)$, $ (1/3,2/9),$ and $(1/3,1/5)$,
respectively.
\end{enumerate}
\end{ex}

\begin{lm}

If $W$ is non-degenerate, then the group $$H:=\{(\alpha_1, \dots,
\alpha_t) \in (\mathbb{C}^*)^t| W (\alpha_1 x_1, \dots, \alpha_t
x_t)=W(x_1, \dots, x_t)\}$$ of diagonal symmetries of $W$ is finite.
In particular, we have $$H \subseteq \mu_{d/k_1}\times \dots \times
\mu_{d/k_t} \cong k_1 \mathbb{Z}/d \times \dots \times k_t
\mathbb{Z}/d,$$ where $\mu_l$ is the group of $l$th roots of unity.
     \end{lm}

\begin{proof}
First write $W = \sum^s_{j=1}W_j$ with $W_j = c_j \prod
x_l^{b_{l,j}}$ and with $c_j \neq 0$.  The uniqueness of the
fractional degrees is equivalent to saying that the matrix $B =
(b_{l,j})$ has rank $t$.  We may as well assume that $B$ is
invertible.  Now write $h = (h_1, \dots, h_t) \in H$, as $h_j =
exp(u_j + v_j i)$ for $u_j\in \R$ is uniquely determined and $v_j
\in \R$ is determined up to integral multiple of $2\pi i$. The
equations $W(h_1x_1, \dots, h_t x_t) = W(x_1, \dots, x_t)$ can now
be written as $B (u+vi) \equiv 0 \pmod{2\pi i \Z}$. The
invertibility of $B$ shows that $u_l = 0$ for all $l$---thus H lives
in $U(1)^t$, and a straightforward argument shows that the number of
solutions $\pmod{2 \pi i \Z}$ to the equation $B (vi) \equiv 0
\pmod{2\pi i\Z}$ is also finite.
\end{proof}

\subsection*{$W$-spin structures on smooth orbicurves}

Let $(\tilde{\Sigma},\z,\m)$ be a smooth orbicurve (or orbifold
Riemann surface) as defined in \cite{S,T}, i.e.,
$(\tilde{\Sigma},\z,\m)$ is a Riemann surface $\Sigma$ with marked
points $\z=\{z_i\}$ having orbifold structure near each marked point
$z_i$ given by a faithful action of ${\Z/m_i}$. In other words, a
neighborhood of each marked point is uniformized by the branched
covering map $z\to z^{m_i}$.  Let $\varrho:\tilde{\Sigma}\to \Sigma$
be the natural projection to the coarse Riemann surface $\Sigma$.

A line bundle $L$ on $\Sigma$ can be uniquely lifted to an orbifold
line bundle on $\tilde{\Sigma}$. We denote the lifted bundle by the
same $L$.

\begin{df}\label{df:log}  Let $K$ be the canonical bundle of $\Sigma$, and
let $$K_{log} := K \otimes  \O(z_1) \otimes \dots \otimes \O(z_k)$$
be the \emph{log-canonical bundle}, that is to say, the bundle whose
holomorphic sections are $1-$ forms holomorphic away from the
special points $\{Z_i\}$ and with at worst simple poles at the
$Z_i$. $K_{log}$ can be thought of as the canonical bundle of the
punctured Riemann surface $\Sigma-\{z_1, \dots, z_k\}$.  Suppose
that $L_1, \dots, L_t$ are orbifold line bundles on $\tilde{\Sigma}$
with isomorphisms $\varphi_j:W_j(L_1, \dots, L_t) \irightarrow
K_{\log}$, where by $W_j(L_1, \dots, L_t)$ we mean the $j$th
monomial
  of $W$ in $L_i$, $$W_j(L_1, \dots, L_t)=L^{\otimes b_{1j}}_1 \otimes \dots,
\otimes L^{\otimes b_{t j}}_t,$$ and where $K_{log}$ is identified
with its pull-back to $\tilde{\Sigma}$. The tuple \linebreak $(L_1,
\dots, L_t, \varphi_1, \dots, \varphi_s)$ is called a \emph{$W$-spin
structure}.
\end{df}

\begin{df}\label{df:RNSLB}
Suppose that the chart of $\tilde{\Sigma}$ at an orbifold point
$z_i$ is $D/(\Z/m)$ with action $e^{\frac{2\pi
i}{m}}(z)=e^{\frac{2\pi i}{m}}z$.  Suppose that the local
trivialization of an orbifold line bundle $L$ is $(D\times
\C)/(\Z/m)$ with the action
\begin{equation}\label{eq:zm-action}
e^{\frac{2\pi i}{m}}(z, w)=(e^{\frac{2\pi i}{m}}z, e^{\frac{2\pi i
v}{m}}w). \end{equation} When $v=0$, we say that $L$ is
\emph{Ramond} at $z_i$. When $v>0$, we say $L$ is
\emph{Neveu-Schwarz (NS)} at $z_i$.

A $W$-spin structure $(L_1,\dots, L_t,\varphi_1, \dots, \varphi_s)$
is called \emph{Ramond} at the point $z_i$ if the group element
$h=(\exp(2\pi i v_1/m),\dots, \exp(2\pi i v_t/m))$ defined by the
orbifold action on the line bundles $L_k$ at $z_i$ acts trivially on
all the line bundles occurring in the monomial $W_j$.  In other
words, the $W$-spin structure is Ramond if there is a monomial
$W_j=c_j \prod x_l^{be,j}$ in $W$ such that for every $l$ with
$b_{l,j}>0$ the line bundle $L_l$ is Ramond at $z_i$.
\end{df}

\begin{rem}\label{rem:desingularize}
If  $L$ is an orbifold line bundle on a smooth orbifold Riemann
surface $\tilde{\Sigma}$, then the sheaf of locally invariant
holomorphic sections of $L$ is locally free of rank one, and hence
dual to a unique orbifold line bundle $|L|$ on $\Sigma$. We also
denote $|L|$ by $\varrho_* L$, and it corresponds to the
desingularization of $L$ \cite{CR1}(Prop 4.2.1). It can be
constructed as follows.

We keep the local trivialization at other places and change it at
the orbifold point $z_i$ by a $\Z/m$-equivariant map $\Psi:
(D-\{0\})\times \C \rightarrow (D-\{0\})\times \C$ by
\begin{equation}(z, w)\rightarrow (z^m, z^{-v}w),
\label{eq:desing-triv}
\end{equation}
where $\Z/m$ acts trivially on the second $(D-\{0\})\times \C$.
Then, we extend $L|_{((D-\{0\})\times \C)}$ to a smooth holomorphic
line bundle over $\Sigma$ by the second trivialization. Since $\Z/m$
acts trivially, this gives a line bundle over $\Sigma$, which is
$|L|$.  Note that if $L$ is Ramond at $z_i$, then $|L|=L$ locally.
When $L$ is  Neveu-Schwarz at $z_i$, then $|L|$ will differ from
$L$.
\end{rem}

\begin{ex}
A smooth orbifold Riemann surface $\tilde{\Sigma}=(\Sigma,\z, \m)$
has a natural orbifold canonical bundle $K_{\tilde{\Sigma}}$, namely
its (orbifold) cotangent bundle. The desingularization is related to
the canonical bundle of $\Sigma$ by
$$|K_{\tilde{\Sigma}}|=K_{\Sigma}\otimes_i \O(-(m_i-1) z_i).$$

On the other hand, the desingularization of the log-canonical bundle
of $\tilde{\Sigma}$ is again the log-canonical bundle of $\Sigma$,
since $K_{log}$ is Ramond at every marked point (the orbifold action
on $dz$ and on $z$ is the same, so it is trivial on $dz/z$).
\end{ex}

Next we study the sections. Suppose that $s$ is a section of $|L|$
having local representative $g(u)$.  Then, $(z, z^v g(z^m))$ is a
local section of $L$. Therefore, we obtain a section
$\varrho^*(s)\in \Omega^0(L)$ which equals $s$ away from orbifold
points, under the identification given by
Equation~\ref{eq:desing-triv}. It is clear that if $s$ is
holomorphic, so is $\varrho^*(s)$. If we start from an analytic
section of $L$, we can reverse the above process to obtain a section
of $|L|$. In particular, $L$ and $|L|$ have isomorphic spaces of
holomorphic sections. In the same way, there is a map $\varrho^*:
\Omega^{0,1}(|L|)\rightarrow \Omega^{0,1}(L)$, where
$\Omega^{0,1}(L)$ is the space of orbifold $(0,1)$-forms with values
in $L$. Suppose that $g(u)d\bar{u}$ is a local representative of a
section of $t\in \Omega^{0,1}(|L|)$. Then $\varrho^*(t)$ has a local
representative $z^v g(z^m) m \bar{z}^{m-1} d\bar{z}$.  Moreover,
$\varrho^*$ induces an isomorphism from $H^1(|L|)\rightarrow
H^1(L)$.

Suppose now that $L^r\cong K_{log}$ with the action of $\Z/m$ on $L$
as in Equation~(\ref{eq:zm-action}). Since $K_{log}$ is Ramond at
every marked point, we must have $r v=l m$ for some  $l$. The
integer $l$ is non-zero precisely when $v$ is, and thus $L$ is
Neveu-Schwarz at $z_i$ if and only if $l>0$.  Moreover, we have
$v<m$, so $l<r$, and of course $\frac{v}{m}=\frac{l}{r}$.  Suppose
that $s\in \Omega^0(|L|)$ with local representative $g(u)$. Then,
$\varrho^*(s^r)$ has local representative $z^{r v }
g^r(z^m)=z^{ml}g^r(z^m)=u^lg^r(u)$. Hence, $s^r \in \Omega^0(K_{log}
\otimes \O ((-l_i) z_i)$, and thus when $l_i>0$, or equivalently,
when $L$ is Neveu-Schwarz at every $z_i$, we have $s^r\in
\Omega^0(K)$.

\begin{rem}\label{rem:zeros-of-powers}
More generally, if $L^r \cong K_{log}$ on a smooth orbicurve with
action of the local group on $L$ defined by $l_i$ (as above) at each
marked point, then we have $$(\varrho_* L)^r = |L|^r = K_{log}
\otimes \O((-l_i)z_i)$$ locally, near $z_i$. \end{rem}

\begin{prop}\label{prop:unique-h}
Let $(L_1, \dots, L_t)$ be a $W$-spin structure on a smooth
orbicurve.  Suppose that the local group $G_z$ of $z$ acts on $L_j$
by $exp(2 \pi i/m)(z,w_j)=(exp(2 \pi i/m)z, exp(2 \pi i
v_{j}/m)w_j)$.  There is a unique element $h \in H$ such that $exp
(2 \pi iv_{j}/m)=h_j=exp(2 \pi i a_j(h))$ with $a_j(h) \in [0,1)$
for every $j$.  Moreover, for every monomial $W_i$, we have
$$W_i(|L_1|,\dots, |L_t|)\cong K_{\log}\otimes
\O(-\sum^t_{j=1}b_{i j}a_j(h)z)$$ near the point $z$.  Letting $h_l$
define the action of the local group $G_{z_l}$ near $z_l$, we have
the global isomorphism
\begin{align*}
W_i(|L_1|),\dots,|L_t|)&\cong K_{log}\otimes
\O(-\sum^k_{l=1}\sum^t_{j=1}b_{i j}a_j (h_l)z_l)\\
&\cong K_\Sigma\otimes \O(-\sum^k_{l=1}\sum^t_{j=1}b_{i j}(a_j
(h_l)-q_j)z_l).
\end{align*}
\end{prop}

\begin{proof}
The existence and uniqueness of $h \in H$ is a straightforward
generalization of the argument for $W=W_{A_{r-1}}$, given above. The
rest is an immediate consequence of the description of $h$ as
$h=(exp(2 \pi i a_1(h),\dots,exp(2 \pi i a_t(h))$ and the
description of $|L_j|$ in terms of the action of the local group
$G_z$ given above.
\end{proof}

\subsection*{$W$-spin equations}

For each monomial $W_i$, let $D=-\sum^k_{l=1}\sum^t_{j=1}b_{i j}(a_j
(h_l)-q_j)z_l$ be a divisor, then there is a canonical meromorphic
section $s_0$ with divisor $D$. This section provides the
identification
$$
K_\Sigma\otimes\O(D)\stackrel{s_0^{-1}}{\cong} K_\Sigma(D),
$$
where $K_\Sigma(D)$ is the sheaf of local, possibly meromorphic,
sections of $K_{\Sigma}$ with zeros (or poles) determined  by $D$.
When at least one of the line bundles occurring in the monomial
$W_i$ is Neveu-Schwarz at $z_l$, then $D$ is not effective. So the
local section of $K_{\Sigma}(D)$ has zeros, and hence is a natural
sub-sheaf of $K_{\Sigma}$. In general, however, it is a sub-sheaf of
$K_{log}$. For each marked point, there is a canonical local section
$\frac{dz}{z}$ of $K_{log}$. Using the isomorphism $\varphi_i$,
there is a local section $t_i$ of $L_i$ with the property $W_i(t_1,
\dots, t_k)=\frac{dz}{z}$. The choice of $t_i$ is unique up to the
action of $H$. We choose a metric on $K_{log}$ with the property
$|\frac{dz}{z}|=\frac{1}{|z|}$. It induces a unique metric on $L_i$
with property $|t_i|=|z|^{-q_i}$. Using the correspondence between
$L_i$ and $|L_i|$, it induces a metric on $|L_i|$ with the behavior
$|e_i|=|z|^{a_i(h)-q_i}$ near a marked point, where $e_i$ is the
corresponding local section of $|L_i|$. In particular, it is a
singular metric at any marked point where $L$ is Ramond.

As before, we assume that $W=\sum W_i=\sum_{i}(c_i \prod_l x_l^{b_{i
l}})$. Let $u_j=\tilde{u}_j e_j$, then it is easy to see that
\begin{equation}
\overline{\frac{\partial W}{\partial u_j}}\in
\overline{K}_{log}\otimes \overline{|L_j|}^{-1}.
\end{equation}
The bundle $\overline{|L_j|}^{-1}$ is isomorphic to $ |L_j|$
topologically. But there is no canonical isomorphism. However, we
can choose an isomorphism compatible with the metric. It induces an
isomorphism $I_1:\Omega(\Sigma, \overline{|L_j|}^{-1}\otimes
\Lambda^{0,1})\rightarrow \Omega(\Sigma, |L_j|\otimes
\Lambda^{0,1})$, such that for a section $v=\tilde{v}e'_j$, we have
$$
I_1(\tilde{v}\bar{e}'_j\otimes d\bar{z})=\tilde{v}|e'_j|^2e_j\otimes
d\bar{z},
$$
where $e'_j$ is the holomorphic basis of $|L_j|^{-1}$ such that
$e'_j\cdot e_j=1$.

It is obvious that $I_1$ is the unique metric-preserving isomorphism
between the corresponding two spaces, and it is independent of the
choice of the local charts.

Since $I_1(\overline{\frac{\partial W}{\partial u_j}})\in
\overline{K}_{log}\otimes |L_j|$, the $W$-spin equation is defined
below as the first order system of the sections $u_1,\dots,u_t$:
$$
\bar{\partial}u_j+I_1\left(\overline{\frac{\partial W}{\partial
u_j}}\right)=0, \quad\mbox{ for all } j=1,\dots, t.
$$

\begin{rem}
The desingularization of the orbifold line bundle $L_j$ induces
isomorphisms $\varrho_j:\Omega^0(|L_j|)\to \Omega^0(L_j)$ and
$\varrho^*_j:\Omega^{0,1}(|L_j|)\to \Omega^{0,1}(L_j)$. It is easy
to see that $\varrho^*$ commutes with $\bar{\partial}$ and
$\frac{\partial\overline{W}}{\partial \bar{u}_j}$, hence the above
$W$-spin equations can be regarded also as equations defined on the
orbicurve. However, we will study the $W$-spin equations in the
resolution line bundle $|L_j|$.
\end{rem}

Without loss of generality, we can assume that the surface $\Sigma$
contains $k$ disjoint closed unit disc $B_1(z_l)$ centered at marked
points $z_l, l=1,\cdots, k$. Take a compact set $\Sigma_0\subset
\Sigma\setminus \cup^k_{l=1} B_{e^{-1}}(z_l)$ such that
$\{\Sigma_0,B_1(z_1),\cdots, B_k(z_k)\}$ forms a cover of $\Sigma$.
Let $\varphi_0,\varphi_1,\cdots,\varphi_k$ be the $C^\infty$
partition functions subordinate to the above cover. We define the
weighted $L^p$ and $L^p_1$ norms of the section $u_i=\tu_i e_i$ in
$B_1(z_l)$ as:
\begin{align*}
&||u_i||_{p;B_1}=\left(\int_{B_1(z_l)} |\tu_i|^p |e_i|^p |dz
d\bar{z}|\right)^{1/p}.\\
&||u_i||_{1,p;B_1(z_l)}=\left(\int_{B_1(z_l)} (|\tu_i|^p+|\partial
\tu_i|+|\bar{\partial} \tu_i|^p)|e_i|^p |dz d\bar{z}|\right)^{1/p}.
\end{align*}
On the interior domain $\Sigma_0$, the norm $|e_i|$ of the base
$e_i$ is not singular, we have the usual definition of Sobolev norm
$||u_i||_{W^p_k(\Sigma_0)}$ of $u_i$.

The global $L^p$ and $L^p_1$ norms are defined as:
\begin{align*}
&||u_i||_{p}=||\varphi_0 u_i||_{W^p_0(\Sigma_0)}+\sum_{l=1}^k ||\varphi_l u_i||_{p;B_1(z_l)}.\\
&||u_i||_{1,p}=||\varphi_0 u_i||_{W^p_1(\Sigma_0)}+\sum_{l=1}^k
||\varphi_l u_i||_{1,p;B_1(z_l)}.
\end{align*}

The weighted Sobolev space $L^p_1(\Sigma, |L_j|)$ is defined as the
closure of $C^\infty_0(\Sigma\setminus \{z_1,\cdots, z_k\}, |L_j|)$
under the norm $||\cdot||_{1,p}$ and $L^p(\Sigma,
|L_j|\otimes\Lambda^{0,1})$ is the closure of
$C^\infty_0(\Sigma\setminus \{z_1,\cdots, z_k\},
|L_j|\otimes\Lambda^{0,1})$ under the norm $||\cdot||_{p}$.

\begin{df} Sections $(u_1,\dots,u_t)$ are said to be
{\em regular }solutions of the $W$-spin equations
\begin{equation}\label{spin-equ0}
\bar{\partial}u_j+I_1\left(\frac{\overline{\partial W}}{\partial
\bar{u}_j}\right)=0,
\end{equation}
if, for each $j$, we have $u_j\in L^2_1(\Sigma, |L_j|)$, and
$I_1(\frac{\partial W}{\partial u_j})\in L^2(\Sigma,|L_j|\otimes
\Lambda^{0,1})$, and $(u_1,\dots,u_t)$ satisfy the $W$-spin
equations almost everywhere.

\end{df}

\begin{rm}
   A regular solution is continuous near a Ramond marked point. Hence,
   it takes a value at a Ramond marked point, which we call the \emph{residue}
   of the solution.
\end{rm}

The spin equation $\bar{\partial}u_j+I_1(\frac{\overline{\partial
W}}{\partial \bar{u}_j})=0$ has different properties near the Ramond
marked points and Neveu-Schwarz marked points. \

Let $u_j=\tu_j e_j$ be a local expression in a coordinate near a
marked point $z_l$. The $W$-spin equation becomes
\begin{equation}\label{equ-NS}
\frac{\bar{\partial} \tu_i}{\partial \bar{z}}+\sum_j
\overline{\frac{\partial W_j(\tu_1,\dots,\tu_t)}{\partial
\tu_i}z^{\Sigma_{s=1}^t b_{j s}(a_s(h_l)-q_s)}}|e'_i|^2=0.
\end{equation}

If all the monomials of $W$ are Ramond, then we have

\begin{equation}
\frac{\bar{\partial} \tu}{\partial \bar{z}}+\overline{\frac{\partial
W(\tu_1,\dots,\tu_t)}{\partial \tu_j}\frac{1}{z}}|e'_j|^2=0.
\end{equation}

In polar coordinates, this equation can be rewritten as
\begin{equation}\label{local-equ1}
\frac{1}{2}r(\frac{\partial}{\partial
r}+\sqrt{-1}\frac{1}{r}\frac{\partial}{\partial
\theta})\tu_j+\overline{\frac{\partial
W(\tu_1,\dots,\tu_t)}{\partial \tu_j}}r^{2q_j}=0.
\end{equation}

\section{Key examples}

    If all the marked points are Neveu-Schwarz, an easy lemma by
    Witten shows that the spin-equation has only the zero solution.
    This was the context people worked with on this subject
    for a long time. It was not clear that in the presence of
    a Ramond marked point the spin equation still should have only the
    zero solution. Indeed, for a while this was what people hoped
    for. In this section, we exhibit examples of infinitely many
    regular nonzero solutions of the  $A_{r-1}$-spin equation degenerating to
    a non-regular solution. All the analysis in this paper was
    designed with this example in mind.

    We first start from local solutions of the $A_{r-1}$ spin equation.

\begin{ex}[\bf A local solution of the  $A_{r-1}$-spin equation near Ramond
marked points]\label{exam-3} Near a Ramond marked point, Equation
(\ref{local-equ1}) becomes

$$
\rho\left(\frac{\partial}{\partial
\rho}+\sqrt{-1}\frac{1}{\rho}\frac{\partial}{\partial
\theta}\right)\tu+2r{\bar{\tu}}^{r-1}\rho^{\frac{2}{r}}=0.
$$

If we assume further that $u$ is a real function and depends only on
the radius $\rho$, then we have

$$
\frac{d\tu}{d\rho}=-2r\tu^{r-1}\rho^{\frac{2}{r}-1}.
$$

Now a special local solution is given by

$$
\tu=({r^2(r-2)}\rho^{\frac{2}{r}}+C)^{-\frac{1}{r-2}},
$$
where $C$ is a positive constant.

An easy computation shows that $\tu\in L^p_1$ if and only if
$p<\frac{2}{1-\frac{1}{r}}$.\

\end{ex}

\begin{ex}[\bf Global solutions]\label{exam-glob} Let $({\mathbb C}\mathbb{P}^1, 3, \z)$ be a marked sphere with two
Ramond marked points and one Neveu-Schwarz marked point with trivial
action $e^{2\pi i}(z,w)=(z,w)$.

We shall construct a sequence of global regular solutions of the
$A_{r-1}$-spin equations.

Let ${\mathbb C}\mathbb{P}^1=U_0\cup U_1$, where
$U_0=\{[Z_0,Z_1]|Z_0\neq 0\}$. Let $z=\frac{Z_1}{Z_0}$ be the affine
coordinate in $U_0$,
 so the Fubini-Study metric is given by
$$
\omega=\frac{\sqrt{-1}}{2\pi}\frac{dz\wedge d\bar{z}}{(1+|z|^2)^2}.
$$
So the induced metric on the canonical bundle is given by
$$
|dz|=1+|z|^2.
$$
In $U_0$ the $A_{r-1}$-spin equation is given by
$$
\frac{\bar{\partial}\tilde{u}_0}{\partial
\bar{z}}+\frac{r}{\bar{z}}(1+|z|^2)^{-\frac{2}{r}}|z|^{\frac{2}{r}}\bar{\tilde{u}}_0^{r-1}=0.
$$
If we only consider the real-valued solution, then the above
equation becomes
$$
\frac{d\tilde{u}_0}{d\rho}=-2r\tilde{u}^{r-1}_0\rho^{\frac{2}{r}-1}(1+\rho^2)^{-\frac{2}{r}}.
$$
Therefore
$$
\tilde{u}_0(\rho)=[2r(r-2)\int^\rho_0\tau^{\frac{2}{r}-1}(1+\tau^2)^{-\frac{2}{r}}d\tau+u_0^{-(r-2)}(0)
]^{-\frac{1}{r-2}}.
$$

It is easy to check that $u=\tilde{u}_0(\frac{dz}{z})^{\frac{1}{r}}$
is really a global solution of the $A_{r-1}$-spin equation. Namely
if we represent $u$ by
 $\tilde{u}_1$ in the other local chart $U_1$, then it also satisfies the
$A_{r-1}$-spin equation. Note that
$\tu_0(z)=\tu_1(z)(-1)^{\frac{1}{r}}$ in $U_0\cap U_1$.

One can easily obtain the relation:
$$
R=\tilde{u}^r_1(0)+\tilde{u}_0^r(0)=u_0^r(0)-(\int^\infty_0(\frac{\tau}{1+\tau^2})^{\frac{2}{r}}\frac{1}{\tau}d\tau
+u_0(0)^{-(r-2)})^{-\frac{r}{r-2}}>0.
$$

Thus, if $R\to \infty$(or $u_0(0)\to \infty$), then
$$
\tilde{u}_0(\rho)\rightarrow
[2r(r-2)\int^\rho_0\tau^{\frac{2}{r}-1}(1+\tau^2)^{-\frac{2}{r}}d\tau]^{-\frac{1}{r-2}},
$$
which is not a regular solution of the $r$-spin equation.
\end{ex}

\section{The $\bar{\partial}$ operator in Weighted Sobolev Spaces }

The Fredholm theory of elliptic operators in weighted Sobolev spaces
has been discussed by many authors (see \cite{Lo1}, \cite{Lo2},
\cite{LM} and references there). In this section, first we recollect
the work of Lockhart and Mcowen for general elliptic operators
defined on a noncompact manifolds with finite ends. Then as an
application of their work, we list the corresponding Fredholm
properties for $\bar{\partial}$, while giving some useful estimates.

\subsection{Lockhart-Mcowen theory}

Suppose $X$ is an $n$-dimensional noncompact manifold without
boundary, containing a compact set $X_0$ such that
$$
X\setminus X_0=\{(\omega, \tau): \omega\in \Omega, \tau\in
(0,\infty)\},
$$
where $\Omega$ is a $n-1$-dimensional closed Riemannian manifold
with a smooth measure $d\omega$.

Let $E,F$ be rank-$d$ vector bundles over $X$. Denote by
$C^\infty(E)$ the set of smooth sections and $C^\infty_0(E)$ the set
of smooth sections with compact set. Choose a finite cover
$\{\Omega_1,\cdots,\Omega_N\}$ of coordinate patches of $\Omega$ and
let $X_\nu=\Omega_\nu\times (0,+\infty)$. We can continue to choose
a covering $X_{N+1},\cdots, X_M$ of coordinate patches of $X_0$ such
that $E$ can be trivialized over $X_\nu, \nu=1,\cdots, N,\cdots, M$.
Let $u=(u_1,\cdots, u_d)$ be a trivialization of a section $u$ with
compact support over $X_\nu$, we can define the norm
$$
||u||_{W^p_s(X_\nu)}:=\sum_{|\alpha|\le s}\sum_{l=1}^d ||D^\alpha
u_l||_{W^p_0(X_\nu)},\;\;(D=-i\partial/\partial x)
$$
where we use the measure $d\omega d\tau$ if $\nu=1,\cdots,N$. Let
$\varphi_1,\cdots, \varphi_{N+M}$ be a set of $C^\infty$ partition
functions subordinate to the cover $X_1,\cdots, X_{N+M}$. We define
a norm on $C^\infty_0(E)$ by
$$
||u||_{W^p_s}=\sum^{N+M}_{\nu=1} ||\varphi_\nu u||_{W^p_s(X_\nu)}
$$
and let $W^p_s(E)$ be the closure of $C^\infty_0(E)$ in this norm.
We can add a weight at infinity to generalize this space. Over
$X_\nu,\nu=1,\cdots, N$ we define the weighted norm
$$
||u||_{W^p_{s,\kappa}( X_\nu)}:=\sum_{|\alpha|\le s}\sum_{l=1}^d
||e^{\kappa\tau}D^\alpha u_l||_{W^p_0(X_\nu)}
$$
and replace $W^p_s(E)$ by $W^p_{s,\kappa}(E)$ whose norm is given
below
$$
||u||_{W^p_{s,\kappa}}=\sum^{N+M}_{\nu=N+1} ||\varphi_\nu
u||_{W^p_s(X_\nu)}+\sum^N_{\nu=1}||\varphi_\nu
u||_{W^p_{s,\kappa}(X_\nu)}.
$$

Suppose $A: C^\infty_0(E)\rightarrow C^\infty_0(F)$ is a
differential operator of order $m$ and is translation invariant in
$\tau>0$. If $d=1$, $A$ has the form
$$
A|_{X_\nu}=\sum^m_{q=0} A^{m-q}(\omega, D_\omega)D^q_{\tau},
$$
where $A^{m-q}(\omega, D_\omega)$ is a differential operator of
order $m-q$ in $\omega\in \Omega_\nu$.

If $d>1$, then $A$ is a $d\times d$ matrix of differential operator
of order $m$ and each entry has the above form. Clearly $A$ is a
bounded operator from $W^p_{s+m,\kappa}$ to $W^p_{s,\kappa}$.

If the vector bundles $E$ and $F$ have the decomposition
$$
E=\oplus^J_{j=1} E_j, F=\oplus^I_{i=1}F_i,
$$
we can generalize the definition of weighted Sobolev spaces as
$$
W^p_{s,\kappa}(E)=\oplus^J_{j=1}W^p_{s_j,\kappa}(E),
W^p_{r,\kappa}(F)=\oplus^I_{i=1}W^p_{r_i,\kappa}(F).
$$
Here $s=(s_1,\cdots,s_J), r=(r_1,\cdots,r_I)$ are multiple indices.

The operator $A$ is also decomposed into
$A_{ij}:C^\infty_0(E_j)\rightarrow C^\infty_0(F_i)$ with order
$s_j-r_i$ (if $s_j<r_i$, we let $A_{ij}=0$). Let
$A^0(x,\xi):E_x\rightarrow F_x$ be the principle symbol of $A$. We
say $A$ is elliptic with respect to $(s,r)$ if $det(A^0(x,\xi))\neq
0$ for any nonzero $(x,\xi)$.

\subsubsection*{Spectrum}

Let $\tilde{X}=\Omega\times {\mathcal R}$ be the full cylinder and
let
$$
\tilde{E}=\oplus^J_{j=1}\tilde{E}_j,
\tilde{F}=\oplus^I_{i=1}\tilde{F}_i
$$
are vector bundles over $\tilde{X}$ with the same rank. Suppose
$$
A=A(\omega,D_\omega,D_\tau):C^\infty_0(\tilde{E})\rightarrow
C^\infty_0(\tilde{F})
$$
be a translation invariant elliptic operator of order $(s,r)$. Then
$$
A: \tilde{W}^p_{s,\kappa}(\tilde{E})\rightarrow
\tilde{W}^p_{r,\kappa}(\tilde{F})
$$
is a bounded operator. Here the weights in
$\tilde{W}^p_{s,\kappa}(\tilde{E})$ and
$\tilde{W}^p_{r,\kappa}(\tilde{F})$ are extended over ${\mathcal
R}$.

Ellipticity and analyticity in $\lambda$ can be used (\cite{AN, AV})
to show
\begin{equation}
A(\omega,D_\omega,\lambda): H^p_s(\tilde{E}|\Omega)\rightarrow
H^p_r(F|\Omega)
\end{equation}
is an isomorphism whenever $\lambda\in \C/\mathcal{C}_A$, where
$\mathcal{C}_A$ is the spectrum of $A(\omega,D_\omega,\lambda)$.
Denote its inverse by $R_A(\lambda)$. It is know that there are only
finitely many spectrum points in any complex strip $\{\lambda:
\kappa_1<Im(\lambda)<\kappa_2\}$ and the eigenspace of each spectrum
point is finite dimensional. Denote by $d(\lambda)$ the dimension of
the eigenspace corresponding to the spectrum point $\lambda$ and let
$\D_A:= \{Im(\lambda)\in \R: \lambda\in \mathcal{C}_A\}$.

Take $f\in C^\infty_0(F)$ and consider its Fourier transformation:
$$
\hat{f}(\omega,\lambda)=\int^{+\infty}_{-\infty} \exp[-i\lambda
\tau] f(\omega,\tau)d\tau.
$$

If $\kappa\in \R/\mathcal{D}_A$, then the operator
\begin{equation}
A^{-1}_\kappa f(\omega,\tau)=\frac{1}{2\pi}\int_{Im \lambda=\kappa}
\exp[i\lambda \tau] R_A(\lambda)\hat{f}(\omega,\lambda) d\lambda
\end{equation}
is a bounded operator and is the inverse operator of $A$. We have
the estimate:
\begin{equation}\label{LM-inho-isom}
||u||_{W^p_{s,\kappa}}\le C||Au||_{W^p_{r,\kappa}},
\end{equation}
where $C$ is constant.

\subsubsection*{A priori estimate}

Returning to the vector bundles $E$ and $F$ over $X$, we double
their restrictions to $\Omega\times \R^+$ to define $\tilde{E}$ and
$\tilde{F}$ on $\tilde{X}$. Using the parametrix method we can
obtain the a priori inequality for any $\kappa\in \R$:
\begin{equation}\label{LM-apri-1}
||u||_{W^p_{s,\kappa}}\le
C(||Au||_{W^p_{r,\kappa}}+||u||_{W^p_{s',\kappa}}) \;(s'_j<s_j)
\end{equation}

\subsubsection*{Fredholm theory}

The a priori estimate is not enough to establish the Fredholm
theory, since the embedding $W^p_{s',\kappa}\to W^p_{s,\kappa}$ is
not compact.

Let $X_1=X_0\cup\{(\omega,\tau):\omega\in \Omega, 0<\tau\le 1\}$,
$\varphi_1\in C_0^\infty(X_1)$, with $\varphi_1=1 $ on $X_0$. Let
$\varphi_2=1-\varphi_1$. By (\ref{LM-inho-isom}) and
(\ref{LM-apri-1}), for any $\kappa\in \R\setminus \mathcal{C}_A$ and
$u\in W^p_{s,\kappa}$ we can obtain
\begin{align}\label{LM-Fred-esti}
&||u||_{W^p_{s\kappa}}\le C(||\varphi_2
Au||_{W^p_{r,\kappa}}+||\varphi_1 A u||_{W^p_{r,\kappa}}\nonumber\\&
+||[\varphi_1,A]u||_{W^p_{r,\kappa}}+||[\varphi_2,A]u||_{W^p_{r,\kappa}}+||\varphi_1
u||_{W^p_{r,\kappa}})\nonumber\\
&\le
C(||Au||_{W^p_{r,\kappa}}+||[\varphi_1,A]u||_{W^p_{r,\kappa}}+||[\varphi_2,A]u||_{W^p_{r,\kappa}}+||\varphi_1
u||_{W^p_{r,\kappa}})
\end{align}

The above inequality shows that $A$ has finite dimensional kernel
and closed range. Furthermore, Lockhart and Mcowen (\cite{LM})
proved the following theorem

\begin{thm}\label{LM-Fred-thm} Suppose $A$ is elliptic with respect to $(s,r)$ and is
translation invariant in $\tau>0$. Then we have the following
conclusions:
\begin{itemize}
\item[(1)] There is a discrete set $\D_A\subset \R$ such that the
operator
$$
A: W^p_{s,\kappa}(E)\rightarrow W^p_{r,\kappa}(F)
$$
is Fredholm if and only if $\kappa\in \R\setminus \D_A$.

\item[(2)] For $\kappa_1, \kappa_2\in \R\setminus\D_A$ with
$\kappa_1<\kappa_2$, there is
$$
i_{\kappa_2}(A)-i_{\kappa_1}(A)=N(\kappa_1,\kappa_2),
$$
where $i_{\kappa_j}$ is the Fredholm index of
$A:W^p_{s,\kappa_j}(E)\rightarrow W^p_{r,\kappa_j}(F)$ and
$N(\kappa_1,\kappa_2):=\sum\{d(\lambda):\lambda\in \mathcal{C}_A
\;\text{with}\;\kappa_1<Im(\lambda)<\kappa_2\})$.

\end{itemize}
\end{thm}

In general, $X$ may have multiple cylindrical ends. Assume that
there are $k$ ends, $X(j)=\{(\omega,\tau): \omega\in
\Omega(j),\tau\in \R^+\}, j=1,\cdots, k$, where $\Omega(j)$ is
$n-1$-dimensional closed Riemannian manifold. Then in different ends
we can choose different weights $\kappa$. Let
$\kappa=(\kappa(1),\cdots, \kappa(k))\in \R^k$ be the weight vector,
then in the same way, we can discuss the weighted Sobolev spaces
$W^p_{s,\kappa}(E)$ and $W^p_{r,\kappa}(F)$ and the differential
operators between them. Let $A(j)=A|_{X(j)}$. Similarly we can
define the corresponding quantities $\lambda(j)\in
\mathcal{C}_{A(j)}, \D_{A{j}}, d(\lambda(j))$.

Define $\D_A=\{(\kappa(1),\cdots, \kappa(k)): \;\text{for at least
one }\;j, \kappa(j)=Im(\lambda(j)),\;\text{where}\;\lambda(j)\in
\mathcal{C}_{A(j)}\}$. Define $\kappa_1\le \kappa_2$ is equivalent
to $\kappa_1(j)\le \kappa_2(j)$ for any $j$. Let
$$
N(\kappa_1,\kappa_2):=\sum \{d(\lambda(j)):\lambda(j)\in
\mathcal(C_{A(j)})\;\text{with}\;\kappa_1(j)<Im(\lambda(j))<\kappa_2(j)\}.
$$

Lochhart and Mcowen(\cite{LM}) had the generalization
\begin{crl} If $A$ is elliptic with respect to $(s,r)$ and is
translation invariant in each end $X(j)$. Then
$A:W^p_{s,\kappa}(E)\to W^p_{r,\kappa}(F)$ is Fredholm if and only
if $\kappa\in \R^k\setminus \D_A$. If $\kappa_1,\kappa_2\in
\R^k\setminus \D_A$ and $\kappa_1\le \kappa_2$, then the change of
Freholm index is given by
\begin{equation}\label{inde-incr}
i_{\kappa_2}(A)-i_{\kappa_1}(A)=N(\kappa_1,\kappa_2).
\end{equation}
Note that when $k>1$, $\D_A$ is not a discrete set but the union of
$(k-1)$-dimensional hyperplanes in $\R^k$.
\end{crl}

\subsection{Fredholm theory for $\bar{\partial}$ operator}

In this section, we will prove some a priori estimate of the
$\bar{\partial}$ operator and then apply Lockhart-Mcowen's theory to
the $\bar{\partial}$ operator in weighted Sobolev space to show that
$\bar{\partial}$ is a Fredholm operator under some mild assumptions.
However, we can't use those estimates of the last section directly.
When transformed into cylindrical coordinates $z=e^{-t-i\theta}$,
the norm $||\cdot||_{1,p}$ is not equivalent to the norm
$||\cdot||_{W^p_{1,\kappa}}$. So we have to deduce our required
estimates.

Let $\bar{\partial} u=f$, where $f\in C^\infty(B_1(z_l),|L_j|\otimes
\Lambda^{0,1})$. Choose cylindrical coordinates
$(z=e^{-t-i\theta})$. Let
$$
u=\tu e_j
$$

$$
f=\tf e_j\otimes d\bar{z}=-\tf e^{i\theta-t}e_j\otimes
(dt-\sqrt{-1}d\theta).
$$
Using cylindrical coordinates, the equation becomes
\begin{equation}\label{equ1}
\frac{\partial \tu }{\partial t}+\sqrt{-1}\frac{\partial \tu
}{\partial \theta}=-2\tf e^{i\theta-t}
\end{equation}

Since $||f||_{p}<\infty$, this is equivalent to
$$
\int^\infty_0 \int_0^{2\pi} |\tf|^p e^{-p(a_j(h_l)-q_j)t-2t}dt
d\theta<\infty.
$$
Let $a_{j,l}=a_j(h_l)-q_j+\frac{2}{p}$, then the integral becomes
$$
\int |\tf e^{-a_{j,l}t}|^p <\infty.\
$$
\subsubsection*{A special solution $u_s$}

Extend $\tf$ symmetrically to $(-\infty, \infty)$, and we get an
equation of the form (\ref{equ1}) defined on the whole real line.

The operator
${\bar{\partial}}^{t,\theta}:=\frac{1}{2}(\partial_t+\sqrt{-1}\partial_\theta)$
 is translation invariant, and the spectrum of
the eigenvalue problem  $(i\lambda+\sqrt{-1}\frac{\partial
}{\partial \theta})\varphi=0$ is $i\Z$. If $\kappa \not\in\Z$, then
by the result of last section, we can construct the inverse
$$
\tu_s =({\bar{\partial}}^{t,\theta})^{-1}_\kappa(-2\tf
e^{i\theta-t})
$$
which satisfies the inequality
\begin{equation}\label{hu-Lp}
||\tu_s ||_{W^p_{1,\kappa}}\le C||2\tf e^{-t}||_{W^p_{0,\kappa}}.
\end{equation}
For convenience, we let $u_s=\tu_s e_j=Q_s\circ (f)$.

Now choosing $\kappa-1=-a_{j,l}$ in (\ref{hu-Lp}), we have
$$
\int^\infty_{-\infty} \int_0^{2\pi} |\tu_s e^{(1-a_{j,l})t}|^p dt
d\theta\le C\int^\infty_{-\infty}\int_0^{2\pi} |\tf
e^{-a_{j,l}t}|^p.
$$
So
$$
\int^\infty_{-\infty}\int_0^{2\pi} |\tu_s e^{t}|^p
e^{-p(a_j(h_l)-q_j)t-2t}dt d\theta\le C\int^\infty_0 \int_0^{2\pi}
|\tf|^p e^{-p(a_j(h_l)-q_j)t-2t}dt d\theta,
$$
which induces
\begin{equation}\label{Lp-equi}
\int_{B_1(0)} |\frac{u_s}{z}|^p |dz d\bar{z}|\le C\int_{B_1(0)}
|f|^p|dz d\bar{z}|.
\end{equation}

On the other hand, by (\ref{hu-Lp}) we have the estimate of
derivatives:
\begin{align}
&\int_{B_1(0)}(|\partial_z \tu_s|^p+|\partial_{\bar{z}} \tu_s|^p
)|z|^{p(a_j(h_l)-q_j)}|dzd\bar{z}|\nonumber\\
&=\int_0^\infty \int_0^{2\pi} (|\partial_t
\tu_s|^p+|\partial_{\theta} \tu_s|^p e^{p t}
)e^{-p(a_j(h_l)-q_j)t} e^{-2t}dtd\theta\nonumber\\
&=\int_0^\infty \int_0^{2\pi} (|\partial_t
\tu_s|^p+|\partial_{\theta} \tu_s|^p
)e^{p(1-a_{j,l})t} dtd\theta\nonumber\\
&\le C\int_{B_1(0)} |f|^p|dz d\bar{z}|\label{Lp-1}.
\end{align}

For the estimate of the $L^p$ norm of $u_s$, we have
\begin{equation}\label{Lp_2}
\int |u_s|^p=\int_{S^1\times (0,\infty)}|\tu|^p
e^{-pa_{j,l}t}\\
\le \int_{S^1\times (0,\infty)}|\tu |^p e^{p(1-a_{j,l})t}\le C\int
|f|^p.
\end{equation}

Combining (\ref{Lp-1}) and (\ref{Lp_2}), we obtain
\begin{equation}\label{Lp-markpt}
||u_s||_{L^p_1(B_1)}\le C||f||_{L^p(B_1)}.
\end{equation}\

Now we apply the ordinary Sobolev embedding theorem to the function
$\tu_s r^c$, where $c:=c_{j,l}:=(a_j(h_l)-q_j)$, to get the weighted
Sobolev embedding inequalities for $u_s$.

If $p\le 2$, then for $ 1<q<\frac{2p}{2-p}$,
\begin{align*}
&||\tu_s r^c||_{W^q_0}\le C||\tu_s r^c||_{W^p_1}\\
&\le C\left(\int |\partial_r(\tu_s
r^c)|^p+\frac{1}{r^p}\left|\frac{\partial}{\partial
\theta}(\tu_s r^c)\right|^p+|\tu_s r^c|^p\right)^{\frac{1}{p}}\\
&=C\left(\int (|\partial_r \tu_s|^p+\frac{1}{r^p}|\partial_\theta
\tu_s|^p+\left|\frac{\tu_s}{r}\right|^p+|\tu_s|^p)r^{c p}\right)^{\frac{1}{p}}\\
&\le C\left(\int (|\partial_r \tu_s|^p+\frac{1}{r^p}|\partial_\theta
\tu_s|^p)r^{c p}+|\tu_s|^p r^{c
p}\right)^{\frac{1}{p}}=C||u_s||_{1,p},
\end{align*}
where the third inequality comes from the relation (\ref{Lp-equi}).
Particularly, when $p=2$, we have for any $1<q<\infty$,
$$
||u_s||_q\le C||u_s||_{1,2}.
$$

If $p>2$, by similar argument we have
$$
||u_s||_{C^\alpha}\le C ||u_s||_{1,p},
$$
where $0<\alpha<1-\frac{2}{p}$.\

In summary, we have
\begin{lm}\label{spec-Lp}
If $f\in L^p(B_1(0), |L_j|\otimes \Lambda^{0,1})$ for $p$ satisfying
the condition $a_{j,l}=a_j(h_l)-q_j+2/p\not\in \Z$, then the special
solution $u_s=Q_s\circ f$ satisfies the following estimates:
\begin{enumerate}
\item[(1)] if $1<p<\infty$, then
\begin{equation}
||u_s||_{1,p;B_1(0)}+||\frac{u_s}{z}||_{p;B_1(0)}\le
C||f||_{p;B_1(0)};
\end{equation}
\item[(2)] if $1<p\le 2$, and $1<q<\frac{2p}{2-p}$, then
\begin{equation}
||u_s||_{q;B_1(0)}\le C||u_s||_{1,p;B_1(0)}\le C||f||_{p;B_1(0)};
\end{equation}
\item[(3)] if $p>2$, and $0<\alpha<1-\frac{2}{p}$, then
\begin{equation}
||{\tu}_s r^c||_{C^\alpha(B_1(0))}\le C||u_s||_{1,p;B_1(0)}\le
C||f||_{p;B_1(0)},
\end{equation}
where $c=a_{j,l}-q_j$.
\end{enumerate}
\end{lm}

\subsubsection*{Estimate of the homogeneous solution}

Let $u=\tu e_j$ satisfy $\bar{\partial}u=0$ in $B_1(0)$, so we have
$\bar{\partial}\tu=0$. We have the interior estimate
$$
||\tu||_{W^s_k(B_1(0))}\le C ||\tu||_{W^s_0(B^+_1(0))}, \forall s>0,
$$
where $B_1^+(0)$ is a ball which is a little bit larger than
$B_1(0)$.

Let $p>1$, then for sufficiently small $s>0$ we have
\begin{align*}
||\tu||^s_{W^s_0(B_1^+(0))}&=\int_{B_1^+(0)}|\tu|^s |z|^{s c} |z|^{-s c}\\
&\le \left(\int |\tu|^p |z|^{c p}\right)^{\frac{s}{p}}\left(\int
|z|^{\frac{-scp}{p-s}}\right)^{\frac{p-s}{p}}\\
\le C(\int |\tu|^p |z|^{c p})^{\frac{s}{p}}=C||u||_{p}^s,
\end{align*}

Therefore, for any $k\ge 0$ and small $s>0$, we have
\begin{equation}
||\tu||_{W^s_k(B_1(0))}\le C||u||_{p; B^+_1(0)}.
\end{equation}

Actually, we have a refined inequality
\begin{equation}
||\tu||_{W^s_k(B_1(0))}\le C||u||_{p; B^+_1(0)\setminus
B_{\frac{1}{2}}(0)},
\end{equation}
since $|u|$ is a subharmonic function and we can use the maximum
principle.

By the Sobolev embedding theorem, we have
$$
||\tu||_{C^k(B_1(0))}\le C ||u||_{p; B^+_1(0)\setminus
B_{\frac{1}{2}}(0)},\forall k\ge 0.
$$

Now using the $C^k$ norm estimate, it is easy to obtain the
following lemma.

\begin{lm}\label{homo-Lp}

Let $\bar{\partial}u=0$ and $u\in L^p(B^+_1(0),|L_j|)$ for $p>1$. We
have the estimate:
\begin{enumerate}
\item[(1)] for any $k\ge 0$ and $1<q<\infty$, there exists a $C$
such that
$$
||\tu||_{W^q_k(B_1(0))}\le C ||u||_{p;B_1^+(0)\setminus
B_{\frac{1}{2}}(0)}.
$$
\item[(2)] if $c\ge 0$, then for $1<q<\infty$, there exists a $C$
such that
$$
||u||_{1,q; B_1(0)}\le C ||u||_{p;B_1^+(0)\setminus
B_{\frac{1}{2}}(0)}.
$$
\item[(3)] if $c<0$, then for $ 1<q<\frac{2}{q_j}$, there exists a
$C$ such that the above inequality in (2) holds.

\end{enumerate}
\end{lm}

Combining Lemma \ref{spec-Lp} and \ref{homo-Lp}, we have

\begin{crl}\label{norm-equa} If $c>0$ at $z_l=0$, then for $1<p<2/(1-\bar{\delta}_j)$, where $\bar{\delta}_j
=\min_{l: c_{jl}>0}(c_{jl})$, there is
\begin{equation}
||u||_{1,p;B_1(0)}+||\frac{u}{z}||_{p;B_1(0)}\le
C||u||_{1,p;B_1(0)}.
\end{equation}
\end{crl}

\begin{proof}
For any $u\in L^p_1(B_1(z_l),|L_j|)$, there is a special solution
$u_s$ such that $\bar{\partial}u_s=\bar{\partial}u$ and satisfies
$$
||u_s||_{1,p;B_1(0)}+||\frac{u_s}{z}||_{p;B_1(0)}\le
C||\bar{\partial}u||_{p;B_1(0)}.
$$
On the other hand, we have $\bar{\partial}(u-u_s)=0$. So by Lemma
\ref{homo-Lp},
$$
||\tu-\tu_s||_{C^0}\le C||u-u_s||_{p}.
$$
Therefore when $c>0$ we have
$$
\int_{B_1} \left| \frac{\tu-\tu_s}{z}\right|^p |z|^{pc}\le
C||u-u_s||_p^p \int_{B_1(0)} |z|^{p(-1+c)} \le C ||u||_{1,p}.
$$
This shows that
$$
||\frac{u}{z}||_{p}\le
||\frac{u-u_s}{z}||_{p}+||\frac{u_s}{z}||_p\le C||u||_{1,p}.
$$
\end{proof}

The following lemma is useful in proving inner compactness.

\begin{lm}\label{lm-norm-Ramo} If $c\le 0$, then for $1<p<\infty$ and
any $u=\tu e_j\in L^p_1(B_1(0), |L_j|)$ satisfying $u(0)=0$, there
is
$$
||u||_{1,p;B_1(0)}+||\frac{u}{z}||_{p;B_1(0)}\le
C||u||_{1,p;B_1(0)}.
$$
\end{lm}

\begin{proof} This lemma is an easy consequence of the
following Hardy inequality (\cite{HLP}):
$$
\int^\infty_0 |\frac{f(t)}{t}|^p t^\epsilon dt\le
\left[\frac{p}{\epsilon -p+1} \right]^p \int^\infty_0 |f'(t)|^p
t^\epsilon dt,
$$
for any $f\in C^\infty_0(0,\infty), \lim_{t\to 0} f(t)=0$ and
$\epsilon\neq p-1.$
\end{proof}

\begin{lm}\label{inho+homo-Lp}
Let $\bar{\partial}u=f$ in $B_1^+(0)$, where $u\in
L^p(B_1^+(0),|L_j|)$ and $f\in L^p(B_1^+(0), |L_j|\otimes
\Lambda^{0,1})$. Then $u\in L^p_1(B_1^+(0),|L_j|)$ and the
inequality
\begin{equation}
||u||_{1,p;B_1(0)}\le C\left(||u||_{p, B_1^+(0)\setminus
B_{\frac{1}{2}}(0)}+||f||_{p,B_1^+(0)}\right)
\end{equation}
holds if the following two conditions are satisfied:
\begin{itemize}
\item $a_{j,l}=a_j(h_l)-q_j+2/p\not\in \Z$;

\item either $c\ge 0$, $1<p<\infty$ or $c<0$, $1<p<\frac{2}{q_j}$.
\end{itemize}
\end{lm}

\begin{proof} Under the assumptions on the parameters $c$ and $p$, one
has
\begin{align*}
&||u||_{1,p;B_1(0)}\le ||u-u_s||_{1,p}+||u_s||_{1,p}\le
C(||u-u_s||_{p, B_1^+(0)\setminus B_{\frac{1}{2}}(0)}+||f||_p)\\ &
\le C(||u||_{p, B_1^+(0)\setminus
B_{\frac{1}{2}}(0)}+||u_s||_p+||f||_p) \le C(||u||_{p,
B_1^+(0)\setminus B_{\frac{1}{2}}(0)}+||f||_{p,B_1^+(0)}).
\end{align*}
\end{proof}

Now by the above lemma, it is easy to obtain the following global
estimate.

\begin{lm}\label{lm-glob-Lp} Let $\bar{\partial}u=f$ on $\Sigma$, where $u\in L^p(\Sigma, |L_j|)$ and
$f\in L^p(\Sigma,|L_j|\otimes \Lambda^{0,1})$. Then $u\in L^p_1$,
and the inequality
\begin{equation}\label{glob-Lp}
||u||_{1,p}\le C(||u||_{L^p(\Sigma\setminus \cup^k_{l=1}
B_{\frac{1}{2}}(z_l))}+||\bar{\partial}u||_p)
\end{equation}
holds if the following two conditions are satisfied:
\begin{itemize}
\item if $a_{j,l}=a_j(h_l)-q_j+2/p\not\in\Z$ for any
$l=1,\cdots,k$. \item either $1<p<\infty$ in the case that
$c_{jl}\ge 0$ at all marked points, or $1<p<\frac{2}{q_j}$ if
$c_{jl}<0$ at some marked point.
\end{itemize}
\end{lm}

Now by the previous result, we will show below that
$\bar{\partial}:L^p_1(\Sigma, |L_j|)\rightarrow
L^p(\Sigma,|L_j|\otimes \Lambda^{0,1})$ is a Fredholm operator.

Under the coordinate transformation $z=e^{-t-i\theta}$, the
neighborhood of a marked point $z_l$ can be viewed as a half
infinite cylinder. We can define the weighted Sobolev space
$W^p_{s,\kappa_{j,l}(p)}$ as in last section on $B_1(z_l)\setminus
\{z_l\}=(S^1\times [0,\infty))$, where
$\kappa_{j,l}(p)=-a_{j,l}=-a_j(z_l)+q_j-2/p$. Similarly we can
define the global space $W^p_{s,\kappa(p)}(\Sigma\times \C)$, where
$\kappa(p)=(\kappa_{j,1}(p),\cdots,\kappa_{j,k}(p))$. Also we have
the expression:
$$
\bar{\partial}=-\frac{1}{2}e^{t-i\theta}(\partial_t+\sqrt{-1}\partial_\theta).
=:-e^{t-i\theta}\bar{\partial}^{t,\theta}.
$$
When changed to the cylinder coordinates, the space
$L^p_1(B_1(z_l))$ is equivalent to the space
$\hat{W}^p_{s,1+\kappa_{j,l}(p)}$ whose norm
\begin{align*}
&||\tu||_{\hat{W}^p_{s,1+\kappa_{j,l}(p)}}\\
&=\{\int (|\tu|^p e^{\kappa_{j,l}pt})+(|\partial_t
\tu|^p+|\partial_\theta \tu|^p)e^{(1+\kappa_{j,l})pt} \}^{1/p}.
\end{align*}

Now it is easy to see that the map
$\bar{\partial}:L^p_1(\Sigma)\rightarrow L^p(\Sigma)$ is equivalent
to the composition of the two maps
$$
\hat{W}^p_{1,1+\kappa(p)}\xrightarrow{\bar{\partial}^{t,\theta}}W^p_{0,1+\kappa(p)}
\xrightarrow{-e^{t-i\theta}}W^p_{0,\kappa(p)}.
$$

Since the map $-e^{t-i\theta}\cdot$ is an isomorphism, hence
$\bar{\partial}$ is Fredholm iff $\bar{\partial}^{t,\theta}$ is
Fredholm and $\ind(\bar{\partial})=\ind(\bar{\partial}^{t,\theta}:
\hat{W}^p_{1,1+\kappa(p)}\rightarrow W^p_{0,1+\kappa(p)})$. Lemma
\ref{lm-glob-Lp} shows that $\bar{\partial}^{t,\theta}$ has
finite-dimensional nullspace and closed image. If $\kappa(p)\not\in
\Z^k$, then by Theorem \ref{LM-Fred-thm} $\bar{\partial}^{t,\theta}:
W^p_{1,1+\kappa(p)}\rightarrow W^p_{0,1+\kappa(p)}$ is Fredholm, in
particular the image
$\bar{\partial}^{t,\theta}(W^p_{1,1+\kappa(p)})$ has
finite-dimensional cokernel. Since
$W^p_{1,1+\kappa(p)}=\hat{W}^p_{1,1+\kappa(p)}\cap
W^p_{0,1+\kappa(p)}$,
$\bar{\partial}^{t,\theta}(\hat{W}^p_{1,1+\kappa(p)})$ also has
finite-dimensional cokernel. Therefore we proved that
$\bar{\partial}$ is Fredholm under the assumption for $p$.

\subsubsection*{Boundary value problem and index computation}

To compute the index of $\bar{\partial}$, we need the index gluing
formula and consider the related boundary value problem.

Let $z_l$ be a marked point, and consider the restriction of the
bundle $|L_j|_{B}$ on the disc $B_1(z_l)$. Assume that $
B_1(z_l)\times \C \rightarrow |L_j|_{B}: (z, w)\rightarrow
\Psi_l(z)w$ is a trivialization such that $\Psi_l(e^{i\theta})\R$
forms a totally real bundle on $S^1_l=\partial B_1(z_l)$. Define the
space
$$
L^{p,B}_1(l):=\{ u\in L^p_1(B_1(z_l)): u(e^{i\theta})\in
\Psi_l(e^{i\theta})\R\}.
$$

Under the $(t, \theta)$ coordinates, this space is equivalent to
$\hat{W}^{p,B}_{1,1+\kappa_{j,l}}$ which is the subspace of
$\hat{W}^p_{1,1+\kappa_{j,l}}$ satisfying the boundary value
condition $\tu(e^{-i\theta})\in \Psi_l(e^{i\theta})\R$. Similarly,
one can define the space $W^{p,B}_{1,1+\kappa_{j,l}}$. Also we can
define the space on the interior
$$
W^{p,B}_1(inn):=\{\tu \in W^p_1(\Sigma\setminus \cup_l B_1(z_l));
\tu(e^{i\theta})\in \Psi_l(e^{i\theta})\R \;\text{for }\;
e^{i\theta}\in S^1_l\}.
$$

The above complex Sobolev spaces can be viewed as real Sobolev
spaces of real 2-dimensional vector functions and the Cauchy-Riemann
operator $\bar{\partial}$ becomes a real linear Fredholm operator.
By totally the same way as in the appendix of [MS], one can prove
the index gluing formulas:
\begin{lm}\label{lm-index}
\begin{itemize}
\item[(1)] $\ind(\bar{\partial}^{t,\theta}:
\hat{W}^p_{1,1+\kappa}\rightarrow
W^p_{0,1+\kappa})=\ind(\bar{\partial}^{t,\theta}: W^{p,B}_{1}(inn)
\rightarrow
W^p_{0}(inn))+\sum^k_{l=1}\ind(\hat{W}^{p,B}_{1,1+\kappa_{j,l}}\rightarrow
W^p_{0,1+\kappa_{j,l}})$

\item[(2)]$\ind(\bar{\partial}^{t,\theta}:
W^p_{1,1+\kappa}\rightarrow
W^p_{0,1+\kappa})=\ind(\bar{\partial}^{t,\theta}: W^{p,B}_{1}(inn)
\rightarrow
W^p_{0}(inn))+\sum^k_{l=1}\ind(W^{p,B}_{1,1+\kappa_{j,l}}\rightarrow
W^p_{0,1+\kappa_{j,l}})$
\end{itemize}
\end{lm}

\begin{thm}\label{thm-Fred} If $1<p<\frac{2}{q_j}$ and $a_j(h_l)-q_j+2/p\neq 1,2$ for any $l$,
then $\bar{\partial}:L^p_1(\Sigma,|L_j|)\rightarrow L^p(\Sigma,
|L_j|\otimes \Lambda^{0,1})$ is a Fredholm operator. In particular,
if $2<p<2/(1-\bar{\delta}_j)$ we have the relation
\begin{align*}
&\ind(\bar{\partial}: L^p_1(\Sigma, |L_j|)\rightarrow L^p(\Sigma,
|L_j|\otimes
\Lambda^{0,1}))=\\
&\ind(\bar{\partial}^{t,\theta}: W^p_{1,1+\kappa}\rightarrow
W^p_{0,1+\kappa})+\#\{z_l: c_{jl}<0\}
\end{align*}
and the index is independent of $p$ in the interval
$(2,2/1-\delta_j)$.
\end{thm}

\begin{proof} We have already proved that $\bar{\partial}: L^p_1(\Sigma, |L_j|)\rightarrow L^p(\Sigma,
|L_j|\otimes \Lambda^{0,1}))$ is a Fredholm operator if $1<p<2/q_j$
and $a_j(h_l)-q_j+2/p\neq 1,2$. Since the index of $\bar{\partial}$
is equal to the index of $\bar{\partial}^{t,\theta}:
\hat{W}^p_{1,1+\kappa}\rightarrow W^p_{0,1+\kappa}$, by Lemma
\ref{lm-index} we only need to compare the index of the operators
$\bar{\partial}^{t,\theta}:
\hat{W}^{p,B}_{1,1+\kappa_{j,l}}\rightarrow W^p_{0,1+\kappa_{j,l}}$
and $\bar{\partial}: W^{p,B}_{1,1+\kappa_{j,l}}\rightarrow
W^p_{0,1+\kappa_{j,l}}$ near each marked point $z_l$. Corollary
\ref{norm-equa} shows that if $2<p<2/1-\bar{\delta}_j$ then
$W^p_{1,1+\kappa_{j,l}}=\hat{W}^{p}_{1,1+\kappa_{j,l}}\cap
W^p_{0,1+\kappa_{j,l}}=\hat{W}^p_{1,1+\kappa_{j,l}}$. Therefore near
marked points with $c_{jl}>0$, the two indices are equal. The rest
case is to compare the indices near marked points with $c_{jl}<0$.

If $c_{jl}<0$ at $z_l$, we have $0<1+\kappa_{j,l}<1$ if $2<p<2/q_j$.
If $p>2$, we have the inclusion $W^p_{0,1+\kappa_{j,l}}\subset
W^2_0$ which implies that
$$
\coker(\bar{\partial}^{t,\theta}:
\hat{W}^{p,B}_{1,1+\kappa_{j,l}}\rightarrow
W^p_{0,1+\kappa_{j,l}})=\coker(\bar{\partial}^{t,\theta}:
W^{p,B}_{1,1+\kappa_{j,l}}\rightarrow W^p_{0,1+\kappa_{j,l}}).
$$

On the other hand, since the group action is trivial for the
resolved bundle $|L_j|\rightarrow S^1\times [0,\infty)$, the
localization $\Psi_l$ obtained by resolving operation satisfies
$\Psi_l(e^{i\theta})\R=\R$. Thus if $\tu\in
\ker(\bar{\partial}^{t,\theta}:
\hat{W}^{p,B}_{1,1+\kappa_{j,l}}\rightarrow
W^p_{0,1+\kappa_{j,l}})$, then $\tu|_{S^1\times \{0\}}$ is a real
function. In particular a real number is an element in the kernel.
By Lemma \ref{lm-norm-Ramo}, we know that
$$
\hat{W}^{p,B}_{1,1+\kappa_{j,l}}\cap \{\tu\in C(S^1\times
[0,\infty)):\tu(\infty)=0\}=W^{p,B}_{1,1+\kappa_{j,l}}.
$$
So combining those consideration, we have
$$
\ind(\bar{\partial}^{t,\theta}:
\hat{W}^{p,B}_{1,1+\kappa_{j,l}}\rightarrow
W^p_{0,1+\kappa_{j,l}})=\ind(\bar{\partial}^{t,\theta}:
W^{p,B}_{1,1+\kappa_{j,l}}\rightarrow W^p_{0,1+\kappa_{j,l}})+1.
$$
By Lemma \ref{lm-index}, we obtain the conclusion.
\end{proof}

\begin{rem} In Theorem \ref{thm-Fred} we only proved an index
transformation formula for $2<p<2/(1-\bar{\delta}_j)$ and have not
considered the case for general $p$ and not computed  the concrete
index. One reason is that the moduli problem we consider is based on
orbicurves---not the resolved curves, for which we can't do the
gluing operation. However the analysis on resolved curves is more
understandable, and the result is easily translated into results on
orbicurves. Therefore, we concentrates only on the analysis of
resolved curves.
\end{rem}

\section{Inner compactness of the solution spaces of $W$-spin
equations}\label{inne-comp}

In this section, we will discuss the compactness problem for the
$W$-spin equations. We will prove if $R$, the sum of the residues of
$W(u_1,\dots, u_t)$ at each Ramond marked point, is finite, then the
corresponding solution space is compact, hence the so-called ``inner
compactness'' holds. However, as shown by Example \ref{exam-glob},
if $R$ is infinite, then the space of the regular solutions is not
compact. The singular solutions of the $W$-spin equations should be
added to compactify the solution space.

Above all we prove that the regular solutions of the $W$-spin
equations lie in $L^p_1$ space for some $p>2$.

Denote by $P_i(u)$ the nonlinear term of the $W$-spin equations
(\ref{spin-equ0}). Then $u_i=u_{i,s}+(u_i-u_{i,s})$, where
$u_{i,s}=-Q_s\circ P_i(u)$ is the special solution we constructed
before. We have the estimate by (2) of Lemma \ref{spec-Lp}
\begin{equation}\label{inho-Lp}
||u_{i,s}||_{q;B_1(z_l)}\le C||u_{i,s}||_{1,2;B_1(z_l)}\le
C||P_i(u)||_{2;B_1(z_l)},
\end{equation}
for any $1<q<\infty$. On the other hand, $u_i-u_{i,s}$ is a
meromorphic section with a possible singularity at the marked
points. Since $u_i-u_{i,s}\in L^2_1$, by the restriction of
integrability, $u_i-u_{i,s}$ should be a holomorphic section by (1)
of Lemma \ref{homo-Lp}.\

There are two cases:

\begin{enumerate}

\item[(1)] if $c_{il}\ge 0$, then by (2) of Lemma \ref{homo-Lp},
$u_i-u_{i,s}$ is $L^q$ integrable for any $q$ with $1<q<\infty$.

\item[(2)] if $c_{il}<0$, then by (3) of Lemma \ref{homo-Lp},
$u_i-u_{i,s}$ is $L^q$ integrable for $1<q<\frac{2}{q_i}$.

\end{enumerate}

So, at least $u_i$ is $L^q$ integrable for
$1<q<\frac{2}{q_i},i=1,\dots,t$. Moreover, by Lemma \ref{homo-Lp},
we have
\begin{equation}\label{homo-L2}
||u_i-u_{i,s}||_{q,B_1(z_l)}\le C||u_i-u_{i,s}||_{2;B_1(z_l)}\le
C(||u_i||_{2;B_1(z_l)}+||u_{i,s}||_{2;B_1(z_l)}),
\end{equation}
for $1<q<\frac{2}{q_i},i=1,\dots,t$. The inequalities
(\ref{homo-L2}) and (\ref{inho-Lp}) induce
\begin{equation}\label{Lq-L2}
||u_i||_{q;B_1(z_l)}\le
C(||u_i||_{2;B_1(z_l)}+||P_i(u)||_{2;B_1(z_l)}),
\end{equation}
for $1<q<\frac{2}{q_i},i=1,\dots,t$.

\begin{rem} We can use the global $L^p$-estimate of the
$\bar{\partial}$ operator in weighted Sobolev space to get the
estimate of the $L^p_1$ norm, but the classical Sobolev embedding
theorem can't be used here to get the estimate of the $L^q$ norm.
\end{rem}

We estimate the norm $\left|\left|\frac{\partial W}{\partial
u_i}\right|\right|^p_p$ for some $p>2$. For simplicity, we take a
monomial $W_l$. Since $\Sigma_j b_{l j}q_j=1$, we have $b_{l
j}q_j<1$ for each $j$. Choose $p,\epsilon$ such that
$0<\epsilon<q_i,2\le p$ and $p(1-\epsilon)<2$.

If $b_{li}\neq 1$, we can choose the H\"older index group
$$
\left(\frac{1-\epsilon}{b_{l1}q_1},\dots,\frac{1-\epsilon}{(b_{l
i}-1)q_i},\dots,\frac{1-\epsilon}{b_{l t}q_t}
,\frac{1-\epsilon}{q_i-\epsilon}\right)
$$
for small $\epsilon$ such that each entry greater than $1$. By the
H\"older inequality, we have
\begin{align}
&\left|\left|\frac{\partial W_l}{\partial
u_i}\right|\right|^p_{p;B_1(z_l)}=\int|u_1|^{p b_{l1}}\dots
|u_i|^{(b_{l i}-1)p}\dots
|u_t|^{p b_{l t}}\nonumber\\
&\le \left(\int
|u_1|^{\frac{p(1-\epsilon)}{q_1}}\right)^{\frac{q_1b_{l1}}{1-\epsilon}}\dots
\left(\int|u_i|^{\frac{p(1-\epsilon)}{q_i}}\right)^{\frac{q_i(b_{l
i}-1)}{1-\epsilon}}\dots \left(\int
|u_t|^{\frac{p(1-\epsilon)}{q_t}}\right)^{\frac{q_t
b_{l t}}{1-\epsilon}}|\Sigma|^{\frac{q_i-\epsilon}{1-\epsilon}}\nonumber\\
&\le C\left(||u_i||_2,||P_i(u)||_2\right)\le \infty.
\end{align}

If $b_{li}=1$, we use, instead, the following H\"older index group
to obtain the analogous estimate:
$$
\left(\frac{1-\epsilon}{b_{l1}q_1},\dots,\frac{1-\epsilon}{b_{l
{i-1}}q_{i-1}},\frac{1-\epsilon}{b_{l
{i+1}}q_{i+1}}\dots,\frac{1-\epsilon}{b_{l t}q_t}
,\frac{1-\epsilon}{q_i-\epsilon}\right).
$$

Thus, if we let $\delta=\min\{q_1,\dots,q_t\}$ and choose $2\le
p<\frac{2}{1-\delta}$, then $\frac{\partial W}{\partial u_i}$ is
$L^p$ integrable for any $i$. Combining the interior estimate and
the estimate near the marked points, we obtain

\begin{lm}\label{Lp1-L2}
Suppose $(u_1,\dots,u_t)$ are solutions of the $W$-spin equations
(\ref{spin-equ0}), then $u_i$ is $L^p_1$ integrable and
$\frac{\partial W}{\partial u_i}$ is $L^p$ integrable for $2\le
p<\frac{2}{1-q_i}$, and there is the estimate
$$
||u_i||_{1,p}\le C\left(||u_i||_p+\left|\left|\frac{\partial
W}{\partial u_i}\right|\right|_p\right)\le
C\left(||u_i||_2,\left|\left|\frac{\partial W}{\partial
u_i}\right|\right|_2\right),
$$
where $C(||u_i||_2,||\frac{\partial W}{\partial u_i}||_2)$ is a
constant depending on the norms $||u_i||_2,$ and
$\left|\left|\frac{\partial W}{\partial u_i}\right|\right|_2$.
\end{lm}

\begin{crl} Suppose $u$ is the solution of an $r$-spin equation, then
$u$ is smooth away from the marked points and is $L^p_1$ integrable
for $2\le p<\frac{2}{1-\frac{1}{r}}$.
\end{crl}

By the above lemma and using the classical Sobolev embedding theorem
near the Ramond marked points (since in this case the weighted norm
control the classical norm), we know each $\tu_i$ is continuous at
marked point $z_l$ with $c_{il}<0$ and have the estimate:
\begin{equation}\label{c0-norm-1}
|\tu_i|_{C^0}\le C||\tu_i||_{W^p_1(B_1(z_l))}\le
C||u_i||_{1,p;B_1(z_l)}.
\end{equation}
Therefore, we can give the following definition.

\begin{df} Suppose that $z$ is a Ramond point and $W_j$ is a Ramond monomial in $W$, then for
sections $u_i, i=1,\dots, t$, $W_j(u_1,\dots,u_t)$ lies in the
log-canonical bundle $K_{log}$. Then Res$W_j$ at $z$ is defined as
the coefficient of the base $\frac{dz}{z}$. If locally we have the
representation $u_i=\tu_i e_i$, then
Res$W_j(u_1,\dots,u_t)|_z=W_j(\tu_1(z),\dots, \tu_t(z))$.
\end{df}

\textbf{Further estimate of the $W$-spin equations}\

Consider the following integral
$$
\sum_i \left(\bar{\partial}u_i, I_1\left(\frac{\overline{\partial
W}}{\partial u_i}\right)\right)_{L^2}
$$
over $\Sigma$.

We will show that the Neveu-Schwarz marked points and the Ramond
marked points have different contributions to the integral. For
simplicity, we assume there is only one marked point on a smooth
curve $\Sigma$.

1. Assume at this marked point $z_l=0$ that the monomial $W_j$ is
Ramond, then
\begin{align*}
&\sum_i (\bar{\partial}u_i, I_1(\frac{\overline{\partial
W_j}}{\partial u_i}))_{L^2(\Sigma)}=\lim_{\epsilon\rightarrow 0}
\sum_i (\bar{\partial}u_i, I_1(\frac{\overline{\partial
W_j}}{\partial u_i}))_{L^2(\Sigma\setminus B_\epsilon(0))}\\
&= \lim_{\epsilon\rightarrow 0}\sum_i \int_{\Sigma\setminus
B_\epsilon(0)} (\frac{\bar{\partial}\tu_i}{\partial
\bar{z}}d\bar{z}\otimes e_i, \overline{\frac{\partial
W_j(\tu_1,\dots,\tu_t)}{\partial
\tu_i}\frac{1}{z}}|e'_i|^2 e_i\otimes d\bar{z})\\
&=\lim_{\epsilon\rightarrow 0}\sum_i\int
\frac{\bar{\partial}\tu_i}{\partial \bar{z}}\frac{\partial
W_j(\tu_1,\dots,\tu_t)}{\partial
\tu_i}\frac{1}{z} dz\wedge d\bar{z}\frac{\sqrt{-1}}{2}\;\;(\text{since}\;*(|e'_i|^2 e_i=e'_i))\\
&=\lim_{\epsilon\rightarrow 0}\frac{\sqrt{-1}}{2}\int
\frac{\bar{\partial}}{\partial
\bar{z}}(W_j(\tu_1,\dots,\tu_t))\frac{1}{z}dz\wedge d\bar{z}\\
&=\lim_{\epsilon\to 0}\frac{\sqrt{-1}}{2}\int_{\partial
B_\epsilon(0)} \frac{W_j(\tu_1,\dots,\tu_t)}{z}dz=-\pi
W_j(\tu_1(0),\dots,\tu_t(0)).
\end{align*}

2. Assume at this marked point that $W_j$ is Neveu-Schwarz.
Furthermore without loss of generality, we can assume that for $1\le
i\le  t_l$ the bundles $|L_i|$ are Ramond and for $t_{l}+1\le i\le
t$ the bundles $|L_i|$ are Neveu-Schwarz.

We have
\begin{align*}
&\sum_i (\bar{\partial}u_i, I_1(\frac{\overline{\partial
W_j}}{\partial u_i}))_{L^2}=\sum_i \int
(\frac{\bar{\partial}\tu_i}{\partial \bar{z}}d\bar{z}\otimes e_i,
\overline{\frac{\partial W_j(\tu_1,\dots,\tu_t)}{\partial
\tu_i}z^{\sum_{s=1}^t b_{j s}(a_s(h_0)-q_s)}}|e'_i|^2 e_i\otimes d\bar{z})\\
&=\lim_{z\to 0}-\pi W_j(\tu_1(z),\dots,\tu_t(z))z^{\sum_{s=1}^t
b_{j s}(a_s(h_0))}.\\
\end{align*}
If $1\le i\le t_l$, then $a_i(h_0)=0$. the $C^0$ norm of $\tu_i$ is
controlled by inequality (\ref{c0-norm-1}). If $t_{l}+1\le i\le t$,
$c_{i0}\le 0$ and $a_i(h_0)>0$, we still have the control of the
$C^0$ norm by (\ref{c0-norm-1}). Assume that $t_{l}+1\le i\le t$ and
$c_{i0}>0$. By Corollary \ref{norm-equa}, if $2\le
p<\frac{2}{1-\bar{\delta}_i}$, then
$$
||u_i||_{1,p;B_1(0)}+||\frac{u_i}{z}||_{p,B_1(0)}\le
C\left(||u_i||_2,\left|\left|\frac{\partial W}{\partial
u_i}\right|\right|_2\right).
$$
This is equivalent to
\begin{equation}\label{smoothns}
||\tu_i r^{c_{i0}}||_{W^p_1(B_1(0))}\le
C\left(||u_i||_2,\left|\left|\frac{\partial W}{\partial
u_i}\right|\right|_2\right),
\end{equation}
for $c_{i0}>0$.

By the Sobolev embedding inequality, we have
$$
|\tu_i(z) r^{c_{i0}}|_{C^0(B_1(0))}\le
C\left(||u_i||_2,\left|\left|\frac{\partial W}{\partial
u_i}\right|\right|_2\right).
$$
Therefore,
\begin{align*}
&|W_j(\tu_1(z),\dots,\tu_t(z))z^{\sum_{s=1}^t b_{j
s}(a_s(h_0))}|\\
&\le C\left(||u_i||_2,\left|\left|\frac{\partial W}{\partial
u_i}\right|\right|_2\right) r^{\min_{t_{l}+1\le i\le t} \{a_i(h_0),
q_i\}}
\end{align*}
So we have
$$
\sum_i \left(\bar{\partial}u_i, I_1\left(\frac{\overline{\partial
W_j}}{\partial u_i}\right)\right)_{L^2}=0.
$$

If $\Sigma$ is a nodal curve, then by a similar argument, we can
prove that the nodal points make no contribution to the integral.

In general, if there are multiple marked points, one has
\begin{align*}
&\sum_i \left(\bar{\partial}u_i, I_1\left(\frac{\overline{\partial
W}}{\partial u_i}\right)\right)_{L^2}=\sum_j \sum_i
\left(\bar{\partial}u_i, I_1\left(\frac{\overline{\partial
W_j}}{\partial u_i}\right)\right)_{L^2}\\
&=-\pi \sum_{z_l:Ramond}\sum_{j: W_j is Ramond}
W_j(\tu_1(z_l),\dots, \tu_t(z_l))\\
&=-\pi \sum_{z_l:Ramond}\sum_{j: W_j is Ramond}
\text{Res}W_j(u_1,\dots, u_t)|_{z_l}.
\end{align*}
Let $R:=\sum_{z_l:Ramond}\sum_{j: W_j is Ramond}
\text{Res}W_j(u_1,\dots, u_t)|_{z_l}$, then we have
$$
0=\sum_i \left(
\bar{\partial}u_i,\bar{\partial}u_i+I_1\left(\frac{\overline{\partial
W}}{\partial u_i}\right)\right)_{L^2}=||\bar{\partial}u||^2_2-\pi R.
$$
Therefore, we obtain
\begin{equation}\label{esti-1}
\sum_i||\bar{\partial}u_i||_2^2=\pi R.
\end{equation}
and so
\begin{equation}\label{esti-2}
\sum_i \left|\left|\frac{\partial W}{\partial
u_i}\right|\right|^2_2=\pi R.
\end{equation}

\begin{thm}\label{thm-resi-ener} Suppose that $u_i,i=1,\cdots,t,$ are regular solutions
of the $W$-spin equation. Let $R:=\sum_{z_l:Ramond}\sum_{j: W_j is
Ramond} \text{Res}W_j(u_1,\dots, u_t)|_{z_l}$, then we have
$$
\sum_i||\bar{\partial}u_i||_2^2=\sum_i \left|\left|\frac{\partial
W}{\partial u_i}\right|\right|^2_2=\pi R.
$$
\end{thm}

\begin{crl}[\bf Witten's lemma] Assume that $W$ is non-degenerate. If all the marked points on $\Sigma$
are Neveu-Schwarz points, then the only regular solution of the
$W$-spin equation is the zero solution.
\end{crl}

By inequality (\ref{smoothns}), we obtain the following proposition
about the smoothness of $\tu_i$ at the marked point $z_0$ with
$c_{i0}>0$:

\begin{prop} Assume $z_0$ is a Neveu-Schwarz point of $|L_i|$ with $c_{i0}>0$, then
$\tu_i |z|^{a_i(h_0)-q_i}$ is continuous at $z_0$.
\end{prop}

\textbf{Controlling norms of $u_i$ by $R$}\

Our aim is to control the suitable norms (Sobolev norms or H\"older
continuous norms) of the solutions $u_i$ by $R$, the sum of residues
of $W$ at Ramond marked points.

\begin{thm}\label{thm-new}
Let $W \in \mathbb{C}[x_1, \dots, x_n]$ be a non-degenerate,
quasi-homogeneous polynomial with weights $q_i:=\wt(x_i)<1$ for each
variable $x_i$, with $i=1, \dots,n$. Then for any $n$-tuple $(u_1,
\dots, u_n) \in \mathbb{C}^n$ we have

\[ |u_i| \leq C \left(\sum^n_{j=1}\left|\frac{\partial W}{\partial x_i}(u_1, \dots,
u_n)\right|+1 \right)^{\delta_i},\] where
$\delta_i=\frac{q_i}{\min_j(1-q_j)}$ and the constant $C$ depends
only on $W$. If $q_i\leq 1/2$ for all $ i \in \{1, \dots, n\}$, then
$\delta_i\le 1$ for all $ i \in \{1, \dots, n\}$. If $q_i<1/2$ for
all $ i \in \{1, \dots, n\}$, then $\delta_i<1$ for all $ i \in \{1,
\dots, n\}$.
\end{thm}

An immediate corollary is the following.

\begin{crl}\label{crl-new}  Let $W \in \mathbb{C}[x_1, \dots, x_n]$ be a non-degenerate,
quasi-homogeneous polynomial with weights $q_i:=\wt(x_i)\le 1/2$. If
$(u_1,\dots, u_n)$ are regular solutions of the $W$-spin equations,
then $u_i$ is $L^p_1$ integrable and $\frac{\partial W}{\partial
u_i}$ is $L^p$ integrable for $2\le p<\frac{2}{1-q_i}$, and we have
the estimate
$$
||u_i||_{1,p}\le C\left(\left|\left|\frac{\partial W}{\partial
u_i}\right|\right|_2+1\right),
$$
where $C$ is a constant independent of $u_i,$ for all $
i=1,\dots,t$.
\end{crl}

The proof of the theorem depends primarily on the following.

\begin{lm}\label{lemma-resultant}
For any non-degenerate, quasi-homogeneous polynomial $W \in
\mathbb{C}[x_1, \dots, x_n]$ and any $n$-tuple $s_1, \dots, s_n \in
\mathbb{C}^n$, the values $(u_1, \dots, u_n) \in \mathbb{C}^n$ that
satisfy $$\frac{\partial W}{\partial x_i}(u_1, \dots, u_n) =s_i$$
also satisfy a quasi-homogeneous polynomial $$p_i(x_i) \in
\mathbb{C}[s_1, \dots, s_n][x_i]$$ whose highest degree term in $x_i
$ is constant (that is, independent of $s_1, \dots, s_n$).
\end{lm}

The polynomial $p_i$ corresponds to a sort of ``resultant'' of the
polynomials $f_i:=\frac{\partial W}{\partial x_i}-s_i$.

\begin{proof}[Proof (of Lemma \ref{lemma-resultant}).]
Since $W$ is quasi-homogeneous, the polynomials $f_i:=\frac{\partial
W}{\partial x_i}-s_i \in \mathbb{C}[x_1, \dots, x_n, s_1, \dots,
s_n]$ are also quasi-homogeneous, of total weight $1-q_i$, provided
$s_i$ is assigned weight $1-q_i$ as well.  Let $X$ denote the closed
subvariety of weighted projective space $\mathbb{P}^{2n-1}_{(q_1,
\dots, q_n, 1-q_1, \dots, 1-q_n)}$ defined by the vanishing of all
the $f_i$: $$X=Z(f_1, \dots, f_n) \subseteq \mathbb{P}^{2n-1}_{(q_1,
\dots, q_n, 1-q_1, \dots, 1-q_n)}.$$ Any point of $X$ of the form
$(a_1; \dots; a_i; \dots; a_n; 0;\dots;0)$  corresponds to a
non-trivial solution of $\frac{\partial W}{\partial X_j}=0$; thus
the linear subspace $$E_i:=Z(x_i, s_1, \dots, s_n) =\{(a_1;
\dots;0;\dots a_n;0; \dots;0)\}\subseteq \mathbb{P}^{2n-1}_{(q_1,
\dots, 1-q_n)}$$ does not intersect $X$ when $W$ is non-degenerate.

The projection $\pi_i$ from $E_i$ to the subspace $\{(0;
\dots;x_i;0;\dots;0;s_1;\dots;s_n)\} \cong \mathbb{P}^n_{(q_1,
1-q_1, \dots, 1-q_n)}$ is a proper morphism of projective varieties,
and thus the image $\pi_i(X) \subseteq \mathbb{P}^n_{(q_i, 1-q_1,
\dots, 1-q_n)} $ is a closed subvariety.  Moreover, $\pi_i(X)$ is
not all of $\mathbb{P}^n_{(q_i,1-q_1,\dots, 1-q_n)}$; otherwise the
point $(1;0;\dots;0) \in \mathbb{P}^n_{(q_i,1-q_1,\dots, 1-q_n)}$
would have a point in $X$ lying over it, and that would contradict
the non-degeneracy of $W$.  Consequently, $\pi_i(X)$ lies in a
hypersurface defined by a quasi-homogeneous polynomial $p_i \in
\mathbb{C}[x_i,s_1,\dots, s_n],$ such that $p_i(1;0\dots; 0)$ is not
zero.

In particular, if we write $p_i$ as a polynomial in $x_i$ with
coefficients in $\mathbb{C}[s_1,\dots, s_n]$, then the leading
coefficient is constant---independent of $s_1,\dots, s_n$.
\end{proof}

Writing $p_i$ as a polynomial in $x_i$, we have
$$p_i(x_i) = \sum_{l=0}^{N}c_l(s_1,\dots,s_n) x_i^{N-l}.$$
To bound the size of $x_i$ in terms of the $s_j$, we must calculate
bounds on the degree in $s_j$ of each coefficient $c_l$. For each
$l$, each term of $c_l$ will be of the form $\alpha
s_1^{\sigma_1}\cdots s_n^{\sigma_n}$ for some non-negative integers
$\sigma_j$ and for $\alpha \in \mathbb{C}$.  Since $p_i$ is
quasi-homogeneous, we have
\begin{align*}
N q_i    &= (N-l)q_i + \sum_j\wt(s_j)\sigma_j \\
        &= (N-l)q_i + \sum_j(1-q_j)\sigma_j.\\
\end{align*}
So
$$l q_i \ge \left(\sum_j \sigma_j\right) \min_j\{1-q_j\}.$$
Letting $\delta_i=\frac{q_i}{\min_j(1-q_j)}$, we have
$$l \delta_i\ge \sum_j \sigma_j.$$ This gives
$$|c_l(s_1,\dots,s_n)| \le K(\sum_j |s_j|+1)^{l\delta_i}$$ for some constant $K$, depending on $l$ and $W$, but independent of all $s_j$.

The following lemma is a simple consequence of the Ger\v{s}gorin
disc theorem and is the final tool that we need to bound the roots
of the polynomial $p_i$.
\begin{lm}\label{lm-ggd}
For any polynomial $f(x) = x^N + \sum_{l=1}^N \alpha_l x^{l-1}$, and
any $N$-tuple of positive real numbers $\rho_1,\dots,\rho_N$, let
$D$ be the maximum of $\rho_l/\rho_{l-1} + \rho_l\alpha_l/\rho_N$
for $N\ge l \ge 2$ and $\rho_1\alpha_1/\rho_N$.  Then the roots of
$f$ lie in the circle
$$ \{|z|\le D\}.$$
\end{lm}
\begin{proof}
The lemma follows immediately from applying the Ger\v{s}gorin disc
theorem \cite[Thm 6.1.1]{HJ} to the $N\times N$ matrix $B A_f
B^{-1}$, where $B$ is the diagonal matrix $B := \diag(\rho_1,\dots,
\rho_N)$, and $A_f$ is the companion matrix \cite[Def 3.3.13]{HJ} of
$f$.
\end{proof}

Applying Lemma~\ref{lm-ggd} with $\rho_l = \left(\sum_j
|s_j|+1\right)^{l\delta_i}$ and $f=p_i$ shows that the roots of
$p_i$ are bounded by $$|x_i| \le C \left(\sum_j
|s_j|+1\right)^{\delta_i}$$ for some constant $C$ that depends only
on $W$. This completes the proof of Theorem~\ref{thm-new}.

\begin{rem}
The non-degeneracy of $W$ is essential to the proof of Theorem
\ref{thm-new}. For example, if $W=u^2v^2 + u^4$ then both partial
derivatives are zero for $u=0$ and $v$ arbitrary; thus we cannot
control $v$ by the partial derivatives.
\end{rem}

\begin{proof}[\bf Proof of Theorem \ref{main-thm1}]\

Suppose that $\{u^n_1,\dots,u^n_t\}$ is a sequence of solutions of
the $W$-spin equation
$$
\bar{\partial}u^n_i+I_1\left(\frac{\partial W}{\partial
u_i}(u^n_1,\dots,u^n_t)\right)=0.
$$
Let $u_i^n=\tu^n_i e_i$ in a local coordinate. We will discuss the
compactness of the solutions in two domains. \

(1). Compactness in the interior domain away from the marked
points.\

In this case, the $W$-spin equations have the following form:
$$
\bar{\partial}\tu^n_i+\overline{\frac{\partial W}{\partial
\tu_i}(\tu_1^n,\dots,\tu_t^n)}\phi=0,
$$
where $\phi$ is a $C^\infty$ function. By Corollary \ref{crl-new},
we have
$$
||\tu^n_i||_{W^p_1(inn)}\le C_R.
$$
Here ``inn'' means the inner domain which has a positive distance to
those marked points. Therefore,  by the standard argument of
compactness, there exists a $C^\infty$ function $\tu_i$ and a
subsequence $\tu^n_i$ (same notation as earlier) such that
$$
\tu^n_i\rightarrow \tu_i\;\text{in}\;C^k\;\text{and
ordinary}\;W^p_k\;\text{norms},
$$
for any integer $k\ge 0$. The $\tu_i$ are certainly solutions of the
$W$-spin equations in the interior part.\

(2). Compactness near marked points.\

Let $|L_i|$ is Neveu-Schwarz at $z_l$ with $c_{i l}>0$. Then by
Corollary \ref{norm-equa} and Corollary \ref{crl-new}, we obtain
\begin{equation}\label{comp-homo}
||u_i^n||_{1,p;B_1(0)}+\left|\left|\frac{u_i^n}{z}\right|\right|_{p,B_1(0)}\le
C_R,
\end{equation}
for $2\le p<\frac{2}{1-\bar{\delta}_i}$. This is equivalent to
\begin{equation}
||\tu^n_i r^{c_{i l}}||_{W^p_1(B_1(0))}\le C_R.
\end{equation}
Obviously this also holds for $c_{il}=0$.

Using the ordinary Sobolev compact embedding theorem, there exists a
$\tu_i$ such that $\tu_i r^{c_{i l}}\in C^{\alpha_i}\cap W^q_0$
(ordinary $q$ norm) for $0<\alpha_i<\bar{\delta}_i,1<q<\infty$, and
\begin{equation}
\tu_i^n r^{c_{i l}}\rightarrow \tu_i r^{c_{i l}}
\;\text{in}\;C^{\alpha'_i},
\end{equation}
where $0<\alpha'_i<\alpha_i$.

If $|L_i|$ is Ramond  or Neveu-Schwarz at $z_l$ with $c_{il}<0$, the
inequality (\ref{comp-homo}) is not true for $2\le
p<\frac{2}{1-\bar{\delta}_i}$; so we can't use the same argument in
this case that we did in the Neveu-Schwarz case. Let
$u^n_{i,s}=-Q_s\circ P_i(u^n)$, since $\tu^n_{i,s}(0)=0$, we have
the decomposition
$$
\tu^n_i=\tu^n_i-\tu^n_i(0)+(\tu^n_i-\tu^n_{i,s})(0).
$$
By Lemma \ref{homo-Lp},
\begin{align}
&|\tu^n_i(0)|=|(\tu^n_i-\tu^n_{i,s})(0)|\le
C||u^n_i-u^n_{i,s}||_{2,B_1(0)}\nonumber\\
&\le C(||u^n_i||_{2,B_1(0)}+||u^n_{i,s}||_{2,B_1(0)})\le C_R.
\end{align}

So there exists a constant $A_i$ such that $\tu^n_i(0)\rightarrow
A_i$ (of course, we take the subsequence as usual).

On the other hand, if $2\le p<\frac{2}{1-q_i}$, we have
\begin{align*}
&||u^n_i-u^n_i(0)||_{1,p;B_1(0)}\le
||u^n_i||_{1,p;B_1(0)}+||u^n_i(0)||_{p,B_1(0)}\\
&\le C_R+\left(\int |z|^{-p q_i}\right)^{\frac{1}{p}}|\tu^n_i(0)|\le
C_R.
\end{align*}
We show that the ordinary $W^p_1$ norm of
$(\tu^n_i(z)-\tu^n_i(0))r^{-q_i}$ can be controlled by its weighted
norm $||\cdot||_{1,p}$. Actually by Lemma \ref{lm-norm-Ramo} we have

\begin{align*}
&||\tu^n_i(z)-\tu^n_i(0)r^{-q_i}||_{W^p_1(B_1(0))} \\
&\le C\left(
||u_i^n(z)-u^n_i(0)||_{1,p;B_1(0)}+||\frac{u_i^n(z)-u^n_i(0)}{z}||_{p;B_1(0)}\right)
\\
&\le C||u_i^n(z)-u^n_i(0)||_{1,p;B_1(0)}\le C_R.
\end{align*}

By the Sobolev compact embedding theorem, there exists a subsequence
and a function $\tilde{v}_i$ such that
$$
(\tu^n_i(z)-\tu^n_i(0))r^{-q_i}\rightarrow \tilde{v}_i r^{-q_i}
\;\text{in}\; C^\alpha, \quad\mbox{ for } 0<\alpha<q_i.
$$
In particular, for every $\varepsilon>0$, there exists $N$ such that
for all $ n>N$,
$$
|\tu^n_i-\tu^n_i(0)-\tilde{v}_i|\le \varepsilon r^{q_i}.
$$
Therefore, for all $ z\in B_1(0)$,
$$
|\tu^n_i(z)-A_i-\tilde{v}_i(z)|\le \varepsilon(1+r^{q_i})\le
C\varepsilon,
$$
i.e., $\tu^n_i\rightarrow A_i+\tilde{v}_i:=\tilde{u}_i$ in
$C^0(B_1(0))$.

We need estimate the following term:
\begin{align}
&\left|\left|\frac{\partial W}{\partial
u_i}(u^n_1,\dots,u^n_t)-\frac{\partial W}{\partial
u_i}(u^m_1,\dots,u^m_t)\right|\right|_p\le \sum_j
\left|\left|\frac{\partial W_j}{\partial
u_i}(u^n_1,\dots,u^n_t)-\frac{\partial W_j}{\partial
u_i}(u^m_1,\dots,u^m_t)\right|\right|_p\nonumber\\
&=\sum_j \left|\left|\frac{\partial W_j}{\partial u_i}(\tu^n_1
r^{c_{1l}},\dots,\tu^n_t r^{c_{t l}})-\frac{\partial W_j}{\partial
u_i}(\tu^m_1 r^{c_{1l}},\dots,\tu^m_t r^{c_{t
l}})\right|\right|_{W^p_0}.
\end{align}

Assume that $|L_i|, 1\le i\le t_l$, is line bundle with $c_{il}\le
0$ at $z_l$, where $0\le t_l\le t$, and $|L_i|$ is line bundle with
$c_{il}>0$ for $i\ge t_l+1$. There are two cases which may happen
for each monomial $W_j$:

Case 1. If $i>t_l$, then
\begin{align}
&\left|\left|\frac{\partial W_j}{\partial u_i}(\tu^n_1
r^{c_{1l}},\dots,\tu^n_t r^{c_{t l}})-\frac{\partial W_j}{\partial
u_i}(\tu^m_1 r^{c_{1l}},\dots,\tu^m_t r^{c_{t
l}})\right|\right|_{W^p_0}\nonumber\\
&\le\left|\left|\left(\frac{\partial W_j}{\partial
u_i}(\tu^n_1,\dots,\tu^{n}_{t_l}, \dots,\tu^n_t r^{c_{t
l}})-\frac{\partial W_j}{\partial u_i}(\tu^m_1
,\dots,\tu^{m}_{t_l},\dots,\tu^m_t r^{c_{t
l}})\right)r^{-\sum_{1\le k\le t_l}b_{jk}q_k}\right|\right|_{W^p_0}\nonumber\\
&\le \left|\frac{\partial W_j}{\partial
u_i}(\tu^n_1,\dots,\tu^{n}_{t_l}, \dots,\tu^n_t r^{c_{t
l}})-\frac{\partial W_j}{\partial u_i}(\tu^m_1
,\dots,\tu^{m}_{t_l},\dots,\tu^m_t r^{c_{t l}})\right|_{C^0} \left(
|z|^{(-1+q_i)p}\right)^{\frac{1}{p}}.
\end{align}

Case 2. If $i\le t_l$, then

\begin{align}
&\left|\left|\frac{\partial W_j}{\partial u_i}(\tu^n_1
r^{c_{1l}},\dots,\tu^n_t r^{c_{t l}})-\frac{\partial W_j}{\partial
u_i}(\tu^m_1 r^{c_{1l}},\dots,\tu^m_t r^{c_{t
l}})\right|\right|_{W^p_0}\nonumber\\
&\le\left|\left|\left(\frac{\partial W_j}{\partial
u_i}(\tu^n_1,\dots,\tu^{n}_{t_l}, \dots,\tu^n_t r^{c_{t
l}})-\frac{\partial W_j}{\partial u_i}(\tu^m_1
,\dots,\tu^{m}_{t_l},\dots,\tu^m_t r^{c_{t
l}})\right)r^{-\sum_{1\le k\le t_l}b_{jk}q_k}\right|\right|_{W^p_0}\nonumber\\
&\le \left|\frac{\partial W_j}{\partial
u_i}(\tu^n_1,\dots,\tu^{n}_{t_l}, \dots,\tu^n_t r^{c_{t
l}})-\frac{\partial W_j}{\partial u_i}(\tu^m_1
,\dots,\tu^{m}_{t_l},\dots,\tu^m_t r^{c_{t l}})\right|_{C^0} \left(
|z|^{(-1+\delta)p}\right)^{\frac{1}{p}}.
\end{align}

Therefore, in either cases, we have
\begin{align}\label{comp-esti-w}
&\left|\left|\frac{\partial W}{\partial
u_i}(u^n_1,\dots,u^n_t)-\frac{\partial W}{\partial
u_i}(u^m_1,\dots,u^m_t)\right|\right|_p\le \sum_j
\left|\left|\frac{\partial W_j}{\partial
u_i}(u^n_1,\dots,u^n_t)-\frac{\partial W_j}{\partial
u_i}(u^m_1,\dots,u^m_t)\right|\right|_p\nonumber\\
&\le \sum_j\left|\frac{\partial W_j}{\partial
u_i}(\tu^n_1,\dots,\tu^{n}_{t_l}, \dots,\tu^n_t r^{c_{t
l}})-\frac{\partial W_j}{\partial u_i}(\tu^m_1
,\dots,\tu^{m}_{t_l},\dots,\tu^m_t r^{c_{t l}})\right|_{C^0}
\left( |z|^{(-1+\delta)p}\right)^{\frac{1}{p}}\nonumber\\
&\le C \sum_j\left|\frac{\partial W_j}{\partial
u_i}(\tu^n_1,\dots,\tu^{n}_{t_l}, \dots,\tu^n_t r^{c_{t
l}})-\frac{\partial W_j}{\partial u_i}(\tu^m_1
,\dots,\tu^{m}_{t_l},\dots,\tu^m_t r^{c_{t l}})\right|_{C^0},
\end{align}
where ${\delta}=\min\{q_1,\cdots,q_t\}$ and $2\le p<
\frac{2}{1-{\delta}}$.

By Lemma \ref{inho+homo-Lp} and (\ref{comp-esti-w}), we have for any
$2\le p<\frac{2}{1-\delta}$, $i\ge t_l$:
\begin{align*}
&||u_i^n-u_i^m||_{1,p;B_1(0)}\le
C\left(||u_i^n-u_i^m||_{p;B_1^+(0)}+\left|\left|\frac{\partial
W}{\partial u_i}(u^n_1,\dots,u^n_t)-\frac{\partial W}{\partial
u_i}(u^m_1,\dots,u^m_t)\right|\right|_{p;B_1^+(0)}\right)\\
&=C \left( ||(\tu_i^n-\tu_i^m)r^{c_{i l}}||_{W^p_0(B_1^+(0))}\right.\\
&+\left.\sum_j \left|\frac{\partial W_j}{\partial
u_i}(\tu^n_1,\dots,\tu^{n}_{t_l}, \dots,\tu^n_t r^{c_{t
l}})-\frac{\partial W_j}{\partial u_i}(\tu^m_1
,\dots,\tu^{m}_{t_l},\dots,\tu^m_t r^{c_{t l}})\right|_{C^0}\right).
\end{align*}

By the $C^0$ convergence of $\tilde{u}^n_k$ for $ 1\le k\le t_l$ and
$\tilde{u}^n_k r^{c_{kl}}$ for $t_l+1\le k\le t$, we know that
$\{u^n_i\}$ is a Cauchy sequence in $L^p_1(B_1(0))$, and
\begin{align*}
&u^n_i\rightarrow u_i\;\text{in}\;L^p_1\\
&\frac{\partial W}{\partial u_i}(u^n_1,\dots,u^n_t)\rightarrow
\frac{\partial W}{\partial u_i}(u_1,\dots, u_t)\;\text{in}\;L^p
\end{align*}
for $2\le p<\frac{2}{1-\delta}$ and $i\ge t_l$.

If $1\le i\le t_l$, we have
\begin{align*}
&||u_i^n-u_i^m||_{1,p;B_1(0)}\le
C\left(||u_i^n-u_i^m||_{p;B_1^+(0)}+\left|\left|\frac{\partial
W}{\partial u_i}(u^n_1,\dots,u^n_t)-\frac{\partial W}{\partial
u_i}(u^m_1,\dots,u^m_t)\right|\right|_{p;B_1^+(0)}\right)\\
&=C \left( |(\tu_i^n-\tu_i^m)|_{C^0}\left(\int |z|^{-p q_i}\right)^{\frac{1}{p}}\right.\\
&+\left.\sum_j \left|\frac{\partial W_j}{\partial
u_i}(\tu^n_1,\dots,\tu^{n}_{t_l}, \dots,\tu^n_t r^{c_{t
l}})-\frac{\partial W_j}{\partial u_i}(\tu^m_1
,\dots,\tu^{m}_{t_l},\dots,\tu^m_t r^{c_{t l}})\right|_{C^0}\right).
\end{align*}

This also shows that $\{u^n_i\}$ is a Cauchy sequence in $L^p_1$ and
\begin{align*}
&u^n_i\rightarrow u_i\;\text{in}\;L^p_1\\
&\frac{\partial W}{\partial u_i}(u^n_1,\dots,u^n_t)\rightarrow
\frac{\partial W}{\partial u_i}(u_1,\dots, u_t)\;\text{in}\;L^p
\end{align*}
for $2\le p<\frac{2}{1-\delta}$.

In summary, $(u_1,\dots,u_t)$ is a solution of the $W$-spin
equations in $B_1(0)$.\
\end{proof}

\section{Compactifying the solution space of the $W$-spin
equation}\label{stro-comp}

As shown in Example \ref{exam-glob}, the space of the regular
solutions may not be compact. To compactify this space, we need to
add those solutions having a singularity at Ramond marked points.
First, we give a definition of those solutions.

\begin{df} The sections $(u_1,\dots,u_t)$ are called the \emph{singular
solutions} of the $W$-spin equations if $u_i\in
L^2_{1,loc}(\Sigma\setminus{\text{Ramond marked
points}}),\frac{\partial W}{\partial u_i}\in
L^2_{loc}(\Sigma\setminus{\text{Ramond}}$ marked points), and they
satisfy the $W$-spin equations pointwise away from the Ramond marked
points and are not regular solutions of the $W$-spin equations.
\end{df}

To compactify the solution space in a suitable topology, we have to
consider the asymptotic behavior of the singular solutions near the
Ramond marked points. First, we deduce some basic estimates of the
general $W$-spin equations.

Assume that $0$ is the unique Ramond marked point in $B_2(0)$
 and that $(u_1,\dots,u_t)$ is a singular solution in $B_2(0)-\{0\}$ of the $W$-spin
 equation. Let $u_j=\tu_j e_j$, then locally the $W$-spin equation
 can be written as

\begin{equation}
\frac{\bar{\partial} \tu_i}{\partial \bar{z}}+\sum_j
\overline{\frac{\partial W_j(\tu_1,\dots,\tu_t)}{\partial
\tu_i}z^{\Sigma_{s=1}^t b_{j s}(a_s(h_l)-q_s)}}|e'_i|^2=0.
\end{equation}

If we set $\tu_j=\varphi_j z^{q_j-a_j}$, then an easy computation
shows that the above $W$-spin equation has the following simple
form:
\begin{equation}\label{local-equ2}
\bar{\partial}\varphi_j+\overline{\frac{\partial W}{\partial
\varphi_j}}=0, \quad\mbox{ for all } j=1,\dots,t.
\end{equation}

Note that the $\varphi_j$'s are only locally defined, though their
norms are well-defined in a neighborhood of the origin.

First, by Equation (\ref{local-equ2}), we have the identity

\begin{equation}\label{local-equ3}
\partial_{\bar{z}}W(\varphi_1,\dots,\varphi_t)=-\sum_i\left|\frac{\partial W}{\partial
\varphi_i}\right|^2.
\end{equation}

Since $W$ is a quasi-homogeneous polynomial, it is easy to show

\begin{equation}\label{homo-gene}
\sum_i q_i x_i\partial_{x_i}W=W(x_1,\dots,x_t).
\end{equation}

By Equation (\ref{local-equ2}) and Identity (\ref{homo-gene}), we
have

\begin{equation}
\sum_i q_i \partial_{\bar{z}}\varphi_i\cdot
\bar{\varphi}_i+\overline{W}=0.
\end{equation}

Taking the derivative $\partial_z$ of the above equation, we have

\begin{equation}
\sum_i q_i \partial_z\partial_{\bar{z}}\varphi_i\cdot
\bar{\varphi}_i+\sum_i
q_i|\partial_{\bar{z}}\varphi_i|^2+\overline{\partial_{\bar{z}}W}=0.
\end{equation}

By Equation (\ref{local-equ2}) and (\ref{local-equ3}), we obtain
\begin{equation}\label{eq-secder}
\sum_i
q_i\partial_z\partial_{\bar{z}}\varphi_i\cdot\bar{\varphi}_i=\sum_i(1-q_i)\left|\frac{\partial
W}{\partial \varphi_i}\right|^2.
\end{equation}

Therefore, by Equation (\ref{eq-secder}), we have the important
equation for the norm $N(z):=\sum_i q_i|\varphi_i|^2$:

\begin{equation}
\Delta(\sum_i q_i|\varphi_i|^2)=\sum_i 8(1-q_i)\left|\frac{\partial
W}{\partial \varphi_i}\right|^2+\sum_i 4q_i(|\partial_z
\varphi_i|^2+|\partial_{\bar{z}}\varphi_i|^2).
\end{equation}

This implies the maximum principle (see[CW]).

\begin{lm}\label{lm-maxi}
For any $p>0,0<\theta<1$ and any $R>0$ such that $B_R(z)\in
B_2(0)-\{0\}$, $z\in B_2(0)-\{0\}$, we have

\begin{equation}
\sup_{B_{\theta R}(z)}N(z)\le C\left(\frac{1}{|B_R(z)|}\int_{B_R(z)}
N(z)^p\right)^\frac{1}{p}.
\end{equation}
\end{lm}

To get a pointwise estimate of the upper bound for the solutions, we
need a uniform local integral estimate of $N(z)$. At first we will
discuss the $A_{r-1}$ case, since it is simpler than the general
case and we can get better estimates. Subsequently, we will consider
the general $W$-case for $W=W(x_1, \dots, x_t)$  a non-degenerate
quasi-homogeneous polynomial with all weights $\wt(x_i)<1/2$.

\textbf{The $A_{r-1}$ case}\

If $W=x^r$, then the $W$-spin equation (\ref{local-equ2}) becomes
the $A_{r-1}$-spin equation:

\begin{equation}\label{sol-phi}
\varphi_{\bar{z}}+r\bar{\varphi}^{r-1}=0.
\end{equation}

\begin{lm}\label{lemm-Ar}
Let $\varphi$ be a solution of (\ref{sol-phi}) in $B_2(0)-\{0\}$,
then there exists a constant $C_r$ only depending on $r$ such that
for any $z\in B_2(0)-\{0\}$,
\begin{align}
&|\varphi(z)|\le C_r|z|^{-\frac{1}{r-2}}\\
&|D^\alpha\varphi(z)|\le C_r |z|^{-\frac{1}{r-2}-|\alpha|}.
\end{align}
\end{lm}

\begin{proof}

Equation (\ref{local-equ3}) becomes
\begin{equation}
(\varphi^r)_{\bar{z}}+r^2|\varphi|^{2(r-1)}=0.
\end{equation}

Let $\psi^\beta$, for  $\beta>0$, be a cut-off function with support
away from the origin.  We have
\begin{equation*}
\int
(\varphi^r)_{\bar{z}}\psi^\beta+r^2|\varphi|^{2(r-1)}\psi^\beta=0.
\end{equation*}

Integrating by parts and using the H\"older inequality, we have
\begin{align*}
&\int r^2|\varphi|^{2(r-1)}\psi^\beta=\int
\varphi^r\beta\psi^{\beta-1}\psi_{\bar{z}}\\
&\le \int |\varphi|^r\beta\psi^{\beta-1}|\psi_{\bar{z}}|\\
&\le \beta\left(\int
|\varphi|^{2(r-1)}\psi^\beta\right)^{\frac{r}{2(r-1)}}\left(\int
(\psi^{\beta\frac{r-2}{2(r-1)}-1}|\psi_{\bar{z}}|)^{\frac{2(r-1)}{(r-2)}}\right)^{\frac{r-2}{2(r-1)}}.\\
\end{align*}

Thus we have
\begin{equation}\label{sol-phi1}
\int |\varphi|^{2(r-1)}\psi^\beta\le C_r \int
\psi^{\beta-\frac{2(r-1)}{(r-2)}}|\psi_{\bar{z}}|^{\frac{2(r-1)}{r-2}}.
\end{equation}
Now we take $\beta=\frac{2(r-1)}{r-2}$ and choose $\psi$ satisfying
the requirement that $\psi=1$ in $B_{\frac{|z|}{4}}(z)$, vanishing
outside $B_{\frac{|z|}{2}}(z)$ and $|\nabla \psi|\le \frac{4}{|z|}$.
Here $z\neq 0$. Thus we obtain from (\ref{sol-phi1}) the following
estimate,
\begin{equation}\label{2(r-1)norm}
\int_{B_{\frac{|z|}{4}(z)}(z)} |\varphi|^{2(r-1)}\le C_r
\int_{B_{\frac{|z|}{2}}(z)}|z|^{-\frac{2(r-1)}{(r-2)}}=C_r|z|^{-\frac{2}{r-2}}.
\end{equation}

Now the estimate (\ref{2(r-1)norm}) and Lemma \ref{lm-maxi} induce
the required pointwise estimate. The derivative estimate comes from
the scaling invariance of the $A_{r-1}$-spin equation, i.e., for
every $\epsilon >0$ the function
$\varphi_\epsilon(z):=\epsilon^{\frac{1}{r-2}}\varphi(\epsilon z)$
also satisfies the spin equation.
\end{proof}

By the pointwise estimate, we can get a uniform $L^p$-estimate:
\begin{crl}\label{coro-Ar}
Let $r\ge 3$ and $1<p<2(r-2)$. If $\varphi$ is the solution of
(\ref{sol-phi}) in $B_2(0)-\{0\}$, then $\varphi$ is an integrable
function in $B_2(0)$, and furthermore
$$
||\varphi||_{p,B_2(0)}\le C,
$$
where $C$ depends only on $r,p$.
\end{crl}

We can also obtain the Harnack inequality for $|\varphi|$.

\begin{lm}[\bf Harnack inequality]
Let $0\le \theta<1$ be a fixed number, $0<\epsilon<1$. Assume that
$\varphi$ is a solution of the equation (\ref{sol-phi}) in
$B_2(0)-\{0\}$, then
\begin{align*}
\sup_{z\in T(\epsilon(1-\theta),\epsilon)}|\varphi(z)|\le
C(r,\theta)\inf_{z\in T(\epsilon(1-\theta),\epsilon)}|\varphi(z)|,
\end{align*}
where $T(\epsilon(1-\theta),\epsilon)$ is the annulus with radius
between $\epsilon(1-\theta)$ and $\epsilon$, and $C(r,\theta)$ is a
constant only depending on $r$ and $ \theta$.
\end{lm}

\begin{proof}
By equation (\ref{sol-phi}), we have (since $\varphi\neq 0$)
$$
(\log \varphi)_{\bar{z}}=-r\bar{\varphi}^{r-1}\varphi^{-1}.
$$
Let
\begin{align*}
g(z)=-\frac{1}{\pi}\int_{T(1-\theta,1)}\frac{-r\bar{\varphi}^{r-1}\varphi^{-1}(\zeta)dv}{\zeta-z},
\end{align*}
then $g_{\bar{z}}=-r\bar{\varphi}^{r-1}\varphi^{-1}$, for $z\in
T(1-\theta,1)$. Since $|\varphi(z)|\le C_r
(1-\theta)^{-\frac{1}{r-2}}$ for $z\in T(1-\theta,1)$, then
$|-r\bar{\varphi}^{r-1}\varphi^{-1}|(z)\le C_r(1-\theta)^{-1}$.
Hence $|g(z)|\le C(r,\theta)$, and $g$ is a H\"older continuous
function in $T(1-\theta,1)$. Let $\hat{\Psi}=\log\varphi-g$, then
$$
\hat{\Psi}_{\bar{z}}=0.
$$
Since $\hat{\Psi}$ is continuous, $\hat{\Psi}$ is an analytic
function. We have $\varphi=e^g e^{\hat{\Psi}}$. Let
$\Psi=e^{\hat{\Psi}}$, then $\varphi=e^g\Psi$, where $g$ is H\"older
continuous and $\Psi$ is analytic. The following estimate holds:
\begin{equation}
|\Psi(z)|\le |e^{-g}\varphi(z)|\le C(r,\theta)=:e^L.
\end{equation}
Since $g$ is bounded, to prove the Harnack inequality for $\varphi$,
we need only prove the Harnack inequality for $\Psi$. Now
$L-\log|\Psi(z)|$ is a nonnegative and harmonic function, so we have
the gradient estimate:
\begin{align*}
|\nabla (L-\log|\Psi(z)|)|\le C(r,\theta)(L-\log|\Psi(z)|)\le
C(r,\theta)L,
\end{align*}
i.e., $|\nabla\log|\Psi(z)||\le C(r,\theta)$, which implies the
Harnack inequality in $T(1-\theta,1)$,
\begin{equation}
\sup_{T(1-\theta,1)}|\Psi(z)|\le C\inf_{T(1-\theta,1)}|\Psi(z)|.
\end{equation}

To prove the Harnack inequality in the annulus
$T(\epsilon(1-\theta),\epsilon)$, we use the scaling invariance of
Equation (\ref{sol-phi}). Namely, if $\varphi$ is the solution of
(\ref{sol-phi}) in $T(\epsilon(1-\theta),\epsilon)$, then
$\varphi_\epsilon(z):=\epsilon^{\frac{1}{r-2}}\varphi(\epsilon z)$
is the solution of (\ref{sol-phi}) in the annulus $T(1-\theta,1)$.
Thus one can easily get the same conclusion in the annulus
$T(\epsilon(1-\theta),\epsilon)$.
\end{proof}

Now by the maximum principle and the Harnack inequality, one can
easily get a convergence corollary:
\begin{crl}

Let $\varphi$ be a solution of the equation (\ref{sol-phi}) in
$B_2(0)-\{0\}$, then either $\varphi(z)$ is bounded near the origin
or $\lim_{z\to 0} |\varphi(z)|=\infty$.
\end{crl}

\textbf{The General $W$-case}\

Take $\psi$ as the cut-off function as defined in the $A_{r-1}$
case. Multiplying the two sides of Equation (\ref{local-equ3}) by
$\psi^\beta$ for $\beta$ sufficiently large, and doing integration
by parts in $B:=B_{\frac{|z|}{4}}(z)$, we have

\begin{equation}
\int_B W \beta \psi^\beta\partial_{\bar{z}}\psi=\sum_i \int_B
\left|\frac{\partial W}{\partial \varphi_i}\right|^2\psi^\beta,
\end{equation}

or

\begin{equation}\label{ineq-wint}
\sum_i \int_B\left|\frac{\partial W}{\partial
\varphi_i}\right|^2\psi^\beta=\int \sum_i \frac{\partial W}{\partial
\varphi_i}\varphi_i q_i\beta \psi^{\beta-1}\partial_{\bar{z}}\psi.
\end{equation}

Let $\delta_i, i=1, \dots, t$ be the indices from Theorem
\ref{thm-new}, then we know that $\delta_i=q_i/\min_j\{1-q_j\}$.
Define $\delta_0=\max_i\{\delta_i\}$. Since we assume that
$q_i=\wt(x_i)<1/2$, we have $\delta_0<1$. By Theorem \ref{thm-new}
and (\ref{ineq-wint}), we have

\begin{align*}
\sum_i \int_B |\frac{\partial W}{\partial \varphi_i}|^2\psi^\beta\le
& C\int_B \left(\sum_i \left|\frac{\partial W}{\partial
\varphi_i}\right|+1\right)^{1+\delta_0}\psi^{\beta-1}|\partial_{\bar{z}}\psi|\\
\le &C\int_B \sum_i \left|\frac{\partial W}{\partial
\varphi_i}\right|^{1+\delta_0}\psi^{\beta-1}|\partial_{\bar{z}}\psi|
+\int_B |\partial_{\bar{z}}\psi|\\
\le &\varepsilon \sum_i \int_B \left|\frac{\partial W}{\partial
\varphi_i}\right|^2\psi^\beta+C_\varepsilon \int
\psi^{(\beta-\frac{1+\delta_0}{2}\beta-1)\frac{2}{1-\delta_0}}|\partial_{\bar{z}}\psi|
^{\frac{2}{1-\delta_0}}+C.
\end{align*}

Thus we obtain
\begin{equation}
\sum_i \int_B\left|\frac{\partial W}{\partial
\varphi_i}\right|^2\psi^\beta\le
C|z|^{-\frac{2\delta_0}{1-\delta_0}}+C.
\end{equation}

So if $|z|<r_0$ for some $r_0$ depending only on $W$, we have for
any $z\in B_{r_0}(0)$,

\begin{equation}\label{ineq-wint1}
\sum_i \int_B \left|\frac{\partial W}{\partial
\varphi_i}\right|^2\psi^\beta\le
C|z|^{-\frac{2\delta_0}{1-\delta_0}}.
\end{equation}

By Theorem \ref{thm-new}, Lemma \ref{lm-maxi} and the integral
estimate (\ref{ineq-wint1}), we have
\begin{align}\label{loca-w}
\sum_i q_i|\varphi_i|^2\le&
\left(\frac{1}{|z|^2}\int_{B_{\frac{|z|}{4}}(z)}\sum_i
|\varphi_i|^{\frac{2}{\delta_0}}\right)^{\delta_0}\nonumber\\
\le &C\left(\frac{1}{|z|^2}\int_{B}\sum_i \left|\frac{\partial
W}{\partial \varphi_i}\right|^2\right)^{\delta_0}\nonumber\\
\le & C|z|^{-\frac{2\delta_0}{1-\delta_0}}.
\end{align}
The above inequality implies the following theorem:

\begin{thm}\label{poin-west}
Suppose $W$ to be a non-degenerate quasi-homogeneous polynomial with
all the fractional degrees (or weights) $q_i<1/2, i=1,\dots,t$. Let
$u_i,i=1,\dots,t$ be the solutions of the $W$-spin equations in
$B_1(0)-\{0\}$, then there exist constants $r_0$ and $C$ only
depending on $W$ such that for any $z\in B_{r_0}(z)$ and all $i$,
$$
|u_i|(z)\le C |z|^{-\kappa_i},
$$
where $\kappa_i=\frac{q_i}{1-2q_i}$.
\end{thm}

\begin{proof}
Since $W$ is a quasi-homogeneous polynomial, for any $\lambda\in
{\Bbb R}$ it satisfies the equality
$$
W(\lambda^{k_1}x_1,\dots,\lambda^{k_t}x_t)=\lambda^dW(x_1,\dots,x_t).
$$
Hence $\frac{\partial W}{\partial x_i}$ is also a quasi-homogeneous
polynomial satisfying
\begin{equation}\label{quas-deri}
\frac{\partial W}{\partial
x_i}(\lambda^{k_1}x_1,\dots,\lambda^{k_t}x_t)=\lambda^{d-k_i}\frac{\partial
W}{\partial x_i}(x_1,\dots,x_t).
\end{equation}

Assume that $\varphi_i$'s are solutions of the equation
(\ref{local-equ2}), i.e.,
$$
\bar{\partial}\varphi_i+\overline{\frac{\partial W}{\partial
\varphi_i}}=0, \forall i=1,\dots,t.
$$
Set $\tilde{\varphi}_i:=\lambda^{k_i}\varphi_i(\lambda^{d-2k_i}z)$,
and then we have
\begin{align*}
&\frac{\bar{\partial}\tilde{\varphi}_i}{\partial\bar{z}}=\lambda^{k_i}
\frac{\bar{\partial}\varphi_i}{\partial\bar{z}}|_{\lambda^{d-2k_i}z}\lambda^{d-2k_i}\\
&=-\lambda^{d-k_i}\overline{\frac{\partial W}{\partial
\varphi_i}}(\varphi_1(\lambda^{d-2k_1}z),\dots,
\varphi_t(\lambda^{d-2k_t}z))\\
&=-\overline{\frac{\partial W}{\partial
\varphi_i}}(\lambda^{k_1}\varphi_1(\lambda^{d-2k_1}z),\dots,\lambda^{k_t}\varphi_t
(\lambda^{d-2k_t}z))\\
&=-\overline{\frac{\partial W}{\partial
\varphi_i}}(\tilde{\varphi}_1,\dots,\tilde{\varphi}_t).
\end{align*}
Here, to derive the third equality, we used the relation
(\ref{quas-deri}). The above calculation shows that
$\tilde{\varphi}_i$'s are also solutions of (\ref{local-equ2}). Thus
by the inequality (\ref{loca-w}), there exists a constant $C$ only
depending on $W$ such that for $z\in B_{r_0}(0)$, there is
$$
\tilde{\varphi}_i(\frac{r_0 z}{|z|})\le C.
$$
Setting $\lambda^{d-2k_i}=|z|/r_0$ in the above inequality, we
obtain
$$
|\varphi_i|\le C|z|^{-\kappa_i},
$$
where $\kappa_i=\frac{q_i}{1-2q_i}$. Notice that
$|u_i|=|u_i|_s=|\tilde{u}_i||z|^{-q_i}=|\varphi_i|$, and we are
done.
\end{proof}

Applying Theorem \ref{poin-west} to the $A_{r-1}$ case, we can
recover the first upper-bound estimate in Lemma \ref{lemm-Ar}, since
in this case the fractional degree $q=1/r$ and $\kappa=1/(r-2)$.
Similarly, only by calculating the fractional degrees and
$\kappa_i$, we can obtain the following pointwise estimate of the
$D_{n+1},E_6,E_7,E_8$ cases.

\begin{crl}\label{crl-Dn}
 If $(u_1,u_2)$ are solutions of the $D_{n+1}$-spin
equation
$$\left\{
\begin{array}{l}
\bar{\partial}u_1+I_1(\overline{nu_1^{n-1}+u_2^2})=0\\
\bar{\partial}u_2+I_1(\overline{2u_1u_2})=0
\end{array}\right.
$$
in $B_1(0)-\{0\}$ for $n\ge 3$, then there exist constants $C$ and
$r_0$ depending  only on $n$ such that for $z\in B_{r_0}(0)$, there
is
$$
|u_1|\le C|z|^{-\frac{1}{n-2}},|u_2|\le C|z|^{-\frac{n-1}{2}}.
$$
\end{crl}

\begin{proof} This is because the fractional degrees of the $D_{n+1}$ polynomial are
$(q_1,q_2)=(\frac{1}{n},\frac{n-1}{2n})$ and
$(\kappa_1,\kappa_2)=(\frac{1}{n-2},\frac{n-1}{2})$.
\end{proof}

Since the $E_6,E_8$ cases are the same as the $A_{r-1}$ case, we
will only write down the corollary for the $E_7$ case.

\begin{crl}\label{crl-E7} If $(u_1,u_2)$ are the solutions of the $E_7$-spin
equation
$$\left\{
\begin{array}{l}
\bar{\partial}u_1+I_1(\overline{3u_1^2+u_2^3})=0\\
\bar{\partial}u_2+I_1(\overline{3u_1u_2^2})=0
\end{array}\right.
$$
in $B_1(0)-\{0\}$, then there exist absolute constants $C$ and $r_0$
such that for any $z\in B_{r_0}(0)$, there is
$$
|u_1|\le C|z|^{-1}, |u_2|\le C|z|^{-\frac{2}{5}}.
$$
\end{crl}

\begin{proof} This is because the fractional degrees
$(q_1,q_2)=(\frac{1}{3},\frac{2}{9})$ and
$(\kappa_1,\kappa_2)=(1,\frac{2}{5})$.
\end{proof}

\begin{proof}[\bf Proof of Theorem \ref{main-thm2}] By Theorem \ref{poin-west},
we have a uniform upper bound for the (regular or singular)
solutions on $\partial\Sigma_\epsilon$. Then we can consider the
following integral
$$
\sum_i \left(\bar{\partial}u_i, I_1\left(\frac{\overline{\partial
W}}{\partial \bar{u}_i}\right)\right)_{L^2}
$$
over $\Sigma_\epsilon$ minus some small discs around the
Neveu-Schwarz points. This integral can be reduced to the sum of
some curve integrals over the circles bounding the small discs and
over $\partial \Sigma_\epsilon$. As in the proof of the previous
section, the line integrals over the circles bounding the small
discs will tend to zero as the radius tends to zero. Only the line
integrals over $\partial \Sigma_\epsilon$ contribute. But by Theorem
\ref{poin-west}, we have the uniform pointwise estimate for any
solutions (singular or regular solutions). So we have
$$
\left|\left|\frac{\partial W}{\partial
u_i}\right|\right|_{2,\Sigma_\epsilon}^2\le C_\epsilon,
$$
where $C_\epsilon$ only depends on $\epsilon$. Now the same method
used to prove inner compactness shows that weak compactness also
holds.
\end{proof}\

\begin{proof}[\bf Proof of Theorem \ref{main-thm3}] This follows from
 weak compactness and the fact that all the (regular or
singular) solutions are uniformly integrable in a small neighborhood
of the Ramond marked points as shown in  Theorem \ref{poin-west}.
\end{proof}







\section*{Acknowledgments} 
We would like to thank E.~Witten for sharing many of his
insights into the spin equations on which this paper is based. We
thank Prof.~J\"urgen Jost and Prof.~Xianqing Li-Jost at the
Max-Planck Institute for Mathematics, Leipzig, for their great help.
We also thank Prof.~Weiyue Ding, Prof.~Jiayu Li and Prof.~Guofang
Wang for many useful suggestions, and Prof.~Ralph Kaufmann for
helpful discussions.

\frenchspacing
\bibliographystyle{plain}

\end{document}